\documentclass[12pt]{article}

\usepackage{subcaption}

\usepackage{float}

\usepackage{multirow}

\usepackage{tikz}

\usetikzlibrary{calc,backgrounds,arrows,matrix}

\usepackage{enumerate}

\usepackage{color}
\definecolor{lightblue}{rgb}{0,0.2,0.5}

\usepackage[colorlinks=true, urlcolor=blue,linkcolor=blue, citecolor=lightblue]{hyperref}

\usepackage[round,sort,comma,numbers]{natbib}

\bibpunct{\textcolor{lightblue}{(}}{\textcolor{lightblue}{)}}{,}{a}{}{;}

\usepackage{amssymb,amsmath}

\usepackage{graphicx}

\usepackage{graphics}
\usepackage{color} 

\DeclareMathAlphabet{\eufrak}{U}{}{}{}
\SetMathAlphabet\eufrak{normal}{U}{euf}{m}{n}
\SetMathAlphabet\eufrak{bold}{U}{euf}{b}{n}

\allowdisplaybreaks

 \def\qu{{\mathord{\mathbb Z}}}

 \def\sZZ{{\rm Z\kern-.45em{}Z}}

 \def\sQQ{{\kern 0.27em \vrule height1.45ex width0.03em depth0em
           \kern-0.30em \rm Q}}
 \def\qu{{\mathchoice
         {\sQQ}
         {\sQQ}
   {\kern 0.225em \vrule height1.05ex width0.025em depth0em \kern-0.25em \rm Q}
   {\kern 0.180em \vrule height0.78ex width0.020em depth0em \kern-0.20em \rm Q}
         }}
 \def\sGG{{\kern 0.27em \vrule height1.45ex width0.03em depth0em
           \kern-0.30em \rm G}}
 \def\gg{{\mathchoice
         {\sGG}
         {\sGG}
   {\kern 0.225em \vrule height1.05ex width0.025em depth0em \kern-0.25em \rm G}
   {\kern 0.180em \vrule height0.78ex width0.020em depth0em \kern-0.20em \rm G}
         }}

 \newtheorem{prop}{Proposition}[section]
 \newtheorem{lemma}[prop]{Lemma}

 \newtheorem{theorem}[prop]{Theorem}
 \newtheorem{remark}[prop]{Remark}

\numberwithin{equation}{section}

\newcommand{\re}{\mathrm{e}}

 \newcounter{hyp}
\usepackage{aeguill}

\newenvironment{Proofy}{\removelastskip\par\medskip
\noindent{\em Proof of Theorem} \rm}{\penalty-20\null\hfill$\square$\par\medbreak}

\newenvironment{Proof}{\removelastskip\par\medskip \noindent{\em Proof.} \rm}{\penalty-20\null\hfill$\square$\par\medbreak}

\def\bprf{\begin{Proof}}
\def\nprf{\end{Proof}}
\def\bdes{\begin{description}}
\def\ndes{\end{description}}

\newtheorem{thm}{Theorem}[section]

\def\bdef{\begin{defn}}
\def\ndef{\end{defn}}
\def\bthm{\begin{thm}}
\def\nthm{\end{thm}}
\def\bprop{\begin{prop}}
\def\nprop{\end{prop}}
\def\brmk{\begin{remark}}
\def\nrmk{\end{remark}}
\def\bexa{\begin{exa}}
\def\nexa{\end{exa}}
\def\blem{\begin{lem}}
\def\nlem{\end{lem}}
\def\bcor{\begin{cor}}
\def\ncor{\end{cor}}
\def\bexe{\begin{exe}}
\def\nexe{\end{exe}}

\newcommand{\ee}{\mathbb{E}}
\newcommand{\E}{\mathbb{E}}

\newcommand{\real}{\mathbb{R}}

\def\og{\leavevmode\raise.3ex
     \hbox{$\scriptscriptstyle\langle\!\langle$~}}
\def\fg{\leavevmode\raise.3ex
     \hbox{~$\!\scriptscriptstyle\,\rangle\!\rangle$}~}

\title{\Huge
  Stochastic SIR L\'evy jump model with heavy tailed increments 
}

\author{
 Nicolas Privault\footnote{
\href{mailto:nprivault@ntu.edu.sg}{nprivault@ntu.edu.sg}
 }
 \qquad Liang Wang\footnote{\href{mailto:liang_wang@ntu.edu.sg}{liang\_wang@ntu.edu.sg}
   (corresponding author).
 }
 \\
\small
Division of Mathematical Sciences
\\
\small
School of Physical and Mathematical Sciences
\\
\small
Nanyang Technological University
\\
\small
21 Nanyang Link, Singapore 637371
}

\allowdisplaybreaks

\usepackage{empheq}

\makeatletter
\newcommand*\rel@kern[1]{\kern#1\dimexpr\macc@kerna}
\newcommand*\widebar[1]{%
  \begingroup
  \def\mathaccent##1##2{%
    \rel@kern{0.8}%
    \overline{\rel@kern{-0.8}\macc@nucleus\rel@kern{0.2}}%
    \rel@kern{-0.2}%
  }%
  \macc@depth\@ne
  \let\math@bgroup\@empty \let\math@egroup\macc@set@skewchar
  \mathsurround\z@ \frozen@everymath{\mathgroup\macc@group\relax}%
  \macc@set@skewchar\relax
  \let\mathaccentV\macc@nested@a
  \macc@nested@a\relax111{#1}%
  \endgroup
}
\makeatother

\let\oldcitet=\citet
\let\oldcitep=\citep
\renewcommand{\cite}[1]{\textcolor[rgb]{0,0,1}{\oldcitet{#1}}}
\renewcommand{\citet}[1]{\textcolor[rgb]{0,0,1}{\oldcitet{#1}}}
\renewcommand{\citep}[1]{\textcolor[rgb]{0,0,1}{\oldcitep{#1}}}

\begin{document}

\maketitle

\baselineskip0.6cm

\vspace{-0.6cm}

\begin{abstract}
This paper considers a general stochastic SIR epidemic model driven by a multidimensional L\'evy jump process with heavy tailed increments and possible correlation between noise components. In this framework, we derive new sufficient conditions for disease extinction and persistence in the mean. Our method differs from previous approaches by the use of Kunita's inequality instead of the Burkholder-Davis-Gundy inequality for continuous processes, and allows for the treatment of infinite L\'evy measures by the definition of new threshold values. An SIR model driven by a tempered stable process is presented as an example of application with the ability to model sudden disease outbreak, illustrated by numerical simulations. The results show that persistence and extinction are dependent not only on the variance of the processes increments, but also on the shapes of their distributions. 
\end{abstract}



\baselineskip0.7cm

\parskip-0.1cm


\section{Introduction}
\noindent
In this paper we consider the classical SIR population model perturbed
by random noise.
The model consists in
$S_t+ I_t + R_t$ individuals submitted to a disease,
where $S_t$, $I_t$ and $R_t$ respectively denote
the numbers of susceptible, infected,
and recovered individuals at time $t\in \real_+$,
which are modeled as
\begin{subequations}
\begin{empheq}[left=\empheqlbrace]{align}
  \label{lsir1}
  dS_t&= ( \Lambda -\mu S_t -\beta S_tI_t) dt+S_{t^-} dZ_1 (t),
  \\[3pt]
  \label{lsir2}
  dI_t&= ( \beta S_tI_t-(\mu+\eta+\varepsilon )I_t ) dt+I_{t^-}dZ_2 (t),
  \\[3pt]
  \label{lsir3}
  dR_t&= ( \eta I_t-\mu R_t ) dt+R_{t^-}dZ_3 (t), 
\end{empheq}
\end{subequations}
 where $Z(t) = (Z_1(t),Z_2(t),Z_3(t))$ is a $3$-dimensional
stochastic process modeling the intensity of random perturbations of the system.
Here,
$\Lambda >0$ denotes the population influx into the susceptible
component,
$\beta >0$ reflects the transmission rate from the susceptible group $S_t$ to infected group $I_t$,
$\mu >0$ represents the nature mortality rate of the three compartments $S_t$, $I_t$ and $R_t$,
$\varepsilon >0$ denotes the death rate of infected individuals
induced by the disease, and
$\eta>0$ is the recovery rate of the epidemic. 
\\

The deterministic version of \eqref{lsir1}-\eqref{lsir3}
with $Z(t) = 0$
has been the object of extensive studies, starting with
\cite{kermack} and \cite{may}, 
where the equilibrium of \eqref{lsir1}-\eqref{lsir3}
in the deterministic case has been characterized by
the basic reproduction number
$$
\mathcal{R}_0=\frac{\beta\Lambda}{\mu(\mu+\varepsilon +\eta)}
$$
such that when $\mathcal{R}_0<1$, the system admits
a Globally Asymptotically Stable (GAS) boundary
equilibrium $E_0=({\Lambda}/ {\mu},0,0)$
called the disease-free equilibrium,
whereas
when $\mathcal{R}_0>1$ there exists
a GAS positive equilibrium
  $$
  E^\ast=(S^\ast,I^\ast,R^\ast)=\left(
  \frac{\mu+\varepsilon +\eta}{\beta},\frac{\mu}{\beta}(\mathcal{R}_0-1),\frac{\eta}{\beta}(\mathcal{R}_0-1) \right),
  $$
  which is called the endemic equilibrium. 
In order to model random variations in population numbers, Brownian noise
has been added to the deterministic 
system in e.g. \cite{beddington}, \cite{tornatore}, \cite{mao2011}, 
to better describe the continuous growth of populations in real ecological systems.
\\

L\'evy jump noise has been first incorporated
into the stochastic Lotka-Volterra population model
in \cite{Bao1}, \cite{Bao2}, where uniform $p$th moment bounds and
asymptotic pathwise estimations have been derived.
Driving processes of the form
\begin{equation}
\nonumber 
Z_i(t) = \varrho_i B_i(t) + \int_0^t \int_0^\infty \gamma_i (z) \tilde{N}(ds,dz),
\qquad
i=1,2,3,
\end{equation}
 where $B_1(t)$, $B_2(t)$, $B_3(t)$ are
 independent standard Brownian motions, 
 $\tilde{N}(ds,dz)$ is a compensated Poisson counting process
 with intensity $ds\nu (dz)$ on $\real_+ \times [0,\infty )$ 
and $\nu (dz)$ is a finite L\'evy measure on $[0, \infty)$,
 have been considered 
  in \cite{amllevy}, \cite{wulia} and \cite{xiaobing}.
In this setting, the asymptotic behavior of solutions
of \eqref{lsir1}-\eqref{lsir3}
around the equilibrium of the corresponding deterministic system
has been studied in \cite{amllevy},
and the threshold of this stochastic SIR model
has been investigated in \cite{wulia}.
The asymptotic behavior of the stochastic solution
of an SIQS epidemic model for quarantine modeling with L\'evy jump diffusion term
has been analyzed in \cite{xiaobing}.
\\

Previously used jump models, 
including \cite{amllevy}, \cite{wulia} and \cite{xiaobing}, 
share the property of being
based on a 
Poisson counting process $N(dt,dz)$
 with finite L\'evy measure $\nu(dz)$ on $[0,\infty )$.
  However, this framework excludes important families
  of L\'evy jump processes having an infinite L\'evy measure, 
  as well as flexible correlation between the random noise components
  of the system \eqref{lsir1}-\eqref{lsir3}.
  In particular, the increments of jump-diffusion models with finite L\'evy measures 
  have exponential tails, see e.g. \S4.3 of \cite{cont},
  and they have a limited potential to model 
  extreme events which usually lead to sudden shifts in
population numbers. 
 \\
  
 
  In this paper, we work in the general setting
  of finite or infinite L\'evy measures
  $\nu$ on $\real$,
  which allows us to consider 
  heavy tailed increments having e.g. power law distributions.
  Indeed, empirical data shows that the jump distribution of population dynamics
  under sudden environmental shocks
  such as earthquakes, tsunamis,
  floods, heatwaves and so on,
  can follow power law distributions,
  see e.g. \cite{spl} and references therein.
 \\
 
We consider a $3$-dimensional L\'evy noise
  $Z(t) = (Z_1(t),Z_2(t),Z_3(t))$
  with L\'evy-Khintchine representation
$$
\E \big[
  \re^{iu_1Z_1(t)+iu_2Z_2(t)+iu_3Z_3(t)}
  \big]
= \exp\left(
-\frac{t}{2} \langle u,\varrho u\rangle_{\real^3}
+ t
\int_{\real^3 \setminus\{0\}}
\big(
\re^{i\langle u,\gamma (z)\rangle_{\real^3}}
- i\langle u,\gamma (z)\rangle_{\real^3}
 - 1
 \big)
 \nu (dz)
 \right),
 $$
 $u=(u_1,u_2,u_3)\in \real^3$, $t\in \real_+$,
 where
 $\varrho = (\varrho_{i,j})_{1\leq i,j \leq 3}$
 is a positive definite $3\times 3$
 matrix,
 the functions $\gamma_i:\real^3 \rightarrow \real$, $i=1,2,3$ are
 measurable functions,
 and $\nu (dz)$ is a $\sigma$-finite measure
 of possibly infinite total mass
 on $\real^3\setminus \{0\}$, such that
\begin{equation}
  \label{infinite}
  \int_{\real^3 \setminus\{0\}} \min (
  |\gamma_i (z)|^2 , 1)\nu(dz)<\infty, \qquad i=1,2,3.
\end{equation}
 see e.g. Theorem~1.2.14 in
 \cite{applebk2}.
 In addition, the process $Z(t) = (Z_1(t),Z_2(t),Z_3(t))$ is known to admit
 the representation
\begin{equation}
\label{sdjakl}
Z_i(t) = B^\varrho_i (t)+\int_0^t \int_{\real^3 \setminus\{0\}}\gamma_i(z)\tilde{N}(ds,dz),
\quad
i =1,2,3,
\end{equation}
where $(B^\varrho_1 (t), B^\varrho_2 (t), B^\varrho_3(t))$ is a $3$-dimensional
Gaussian process with independent and stationary increments
and covariance matrix $\varrho = (\varrho_{i,j})_{1\leq i,j \leq 3}$,
and
$\tilde{N}(dt,dz) =N(dt,dz)-\nu(dz)dt$
is the compensated Poisson counting process with L\'evy measure $\nu (dz)$ on
$\real^3\setminus \{0\}$.
All processes are defined on a complete filtered probability space
$(\Omega,\mathcal{F},(\mathcal{F}_t)_{t\geq0},\textbf{P})$,
 $N(dt,dz)$ is
independent of $(B^\varrho_1 (t), B^\varrho_2 (t), B^\varrho_3(t))$,
and the covariances of $(Z_1(t),Z_2(t),Z_3(t))$ are given by
$$
\E [ Z_i(t) Z_j(t) ]
=
\varrho_{i,j} t
+t \int_{\real^3 \setminus\{0\}}\gamma_i(z)\gamma_j(z) \nu (dz ),
\quad
t\in \real_+, 
\quad 
i,j =1,2,3,
$$
which allows for the modeling of random interactions between
the components $(S_t,I_t,R_t)$ of the model.
\\

  In order to investigate the threshold of the stochastic SIR model
  with finite L\'evy measures,
  \cite{wulia} have derived long-term estimates (Lemmas~2.1-2.2 therein)
  which rely on the finiteness of the quantity
    \begin{equation}
    \label{lmbda}
         \int_0^\infty 
 \big( (1+\overline{\gamma}(z))^p-1-\underline{\gamma} (z) \big)\nu(dz), \quad
 p>1,
\end{equation}
 where 
 $$
 \overline{\gamma} ( z) := \max ( \gamma_1(z), \gamma_2(z) , \gamma_3(z) )
 \mbox{ ~and~ }
 \underline{\gamma} ( z) := \min ( \gamma_1(z) , \gamma_2(z) , \gamma_3(z)),
 \quad z\in \real.
 $$
 In our generalized setting
\eqref{sdjakl} under \eqref{infinite}, estimates on solutions 
are obtained by replacing \eqref{lmbda} 
with the expression 
\begin{equation}
\label{lmbda2}
\lambda (p):=
c_p \int_{\real^3 \setminus\{0\}} \overline{\gamma}^2(z) \nu(dz)
+ c_p \int_{\real^3 \setminus\{0\}} \overline{\gamma}^p(z) \nu(dz), \quad p >1.
\end{equation}
 where $c_p:= p (p-1)\max ( 2^{p-3},1 )/2$.
 In addition,
 given the jump stochastic integral process
\begin{equation}
\nonumber 
K_t:= \int_0^t \int_{\real^3 \setminus \{0\}} g_s (z)\ \big(N(ds,dz)-\nu (dz)ds\big),
\quad t\in \real_+,
\end{equation}
of the predictable integrand
$(g_s(y))_{(s,y)\in \real_+ \times \real}$,
 we use Kunita's inequality
\begin{eqnarray}
\label{eq:BDGPoisson2}
\lefteqn{
  \E\left[\sup\limits_{0\leq s\leq t}|K_s|^p\right]
}
\\
\nonumber
& \leq &
   {C}_p
\ee\left[
  \int_0^t\int_{\real^3 \setminus \{0\}} |g_s(z)|^p \ \nu (dz)ds
  \right]
+
{C}_p
\ee\Bigg[ \left(
  \int_0^t \int_{\real^3 \setminus \{0\}} |g_s(z)|^2 \ \nu (dz)ds\right)^{ p/2 }
  \Bigg],
\end{eqnarray}
for all $t\in \real_+$ and $p \geq 2$, where
${C}_p :=
2^{2p-2}
\big(
\sqrt{e^{\log_2p}/p}
+
8 p^{\log_2p}
\big)$,
 see Theorem~2.11 of \cite{kunitalevy},
 Theorem~4.4.23 of \cite{applebk2},
 and Corollary~2.2 in \cite{bretonprivault3}.
 This replaces the
  Burkholder-Davis-Gundy inequality
   for continuous martingales 
   \begin{equation}
     \label{bdg0}
  \E\left[\sup\limits_{0\leq s\leq t}|M_s|^p\right]\leq C_p\E\big[
    \langle M,M \rangle_t^{p/2} \big],
  \qquad p>1,
\end{equation}
   which is used in the proof of Lemmas~2.1 and 2.2 in \cite{wulia}, 
   where $\langle M,M \rangle_t$ is the (predictable)
   quadratic variation of the continuous martingale $(M_t)_{t\in \real_+}$,
   see e.g. Theorem~7.3 in Chapter~1 of \cite{mao2008}.
    Indeed, it is known, see e.g.
Remark~357 in \cite{situ} and \cite{bretonprivault3},
that \eqref{bdg0}
is invalid for martingales with jumps.
\\
 
  As an example, we consider 
  the tempered stable distribution
  introduced in \cite{koponen}
  which belongs to the
    family of self-decomposable distributions,
    see \S~3.15 in \cite{sato}.
 Given $\alpha \in (0,2)$, the tempered $\alpha$-stable L\'evy measure
is defined as 
\begin{equation}\label{measure2}
 \nu (A)=
 \int_0^\infty
 \frac{\re^{-r}}{r^{\alpha +1}}
 \int_{\real^3\setminus \{0\}}
 {\bf 1}_A (rx) R_\alpha (dx) dr,
\end{equation}
 where $R_\alpha (dx)$ is a measure on $\real^3 \setminus \{0\}$ such that
$$
\int_{\real^3 \setminus\{0\}} \min (
\Vert x\Vert^2_{\real^3},
\Vert x \Vert^\alpha_{\real^3}) R_\alpha (dx) <\infty,
$$
see Theorem~2.3 in \cite{rosinski4}.
It has been shown in \cite{sztonik}, Theorem~5, 
    that the increments of tempered stable processes
    can have (heavy) power tails instead of (semi-heavy) exponential tails,
    see also \cite{kuchler}. 
    Taking, e.g.
$$
R_\alpha (dx) = k_- \lambda_-^\alpha
\delta_{(-1/\lambda_-,-1/\lambda_-,-1/\lambda_-)}(dx)
+ k_+ \lambda_+^\alpha \delta_{(1/\lambda_+,1/\lambda_+,1/\lambda_+)}(dx),
$$
where $k_-,k_+,\lambda_- ,\lambda_+ >0$, 
and $\delta_y$ denotes the Dirac measure at $y\in \real^3$,
 the L\'evy measure of the
 $3$-dimensional fully correlated tempered stable process
 is given by 
\begin{eqnarray}
  \label{measure2.0}
  \lefteqn{
    \nu (A)
 =
 \int_{\real^3\setminus \{0\}}
 \int_0^\infty {\bf 1}_A ( r x)
 \frac{\re^{-r}}{r^{\alpha +1}}
 dr R_\alpha (dx)
  }
 \\
 \nonumber
 &  = &
 k_-
  \int_0^\infty {\bf 1}_A
  (-r/\lambda_-,-r/\lambda_-,-r/\lambda_-)
 \frac{\re^{-r}}{r^{\alpha +1}}
 dr
 +
k_+
   \int_0^\infty
       {\bf 1}_A
       (r/\lambda_+ , r/\lambda_+ ,r/\lambda_+ )
 \frac{\re^{-r}}{r^{\alpha +1}} dr, 
\end{eqnarray}
 with $\nu (\real )=+\infty$ for all $\alpha \in (0,2)$.
 We note that 
 \eqref{lmbda} is infinite when 
 $\alpha \in [1,2)$ and $p>1$,
     whereas $\lambda (p)$ given by \eqref{lmbda2}
     remains finite whenever $p > \alpha$. 
     \\

This paper is organised as follows.
After stating preliminary results on the existence and uniqueness
of solutions to \eqref{lsir1}-\eqref{lsir3}
in Proposition~\ref{Theorem 2.1},
in the key Lemmas~\ref{l3.1} and \ref{l3.4} we derive
new solution estimates by respectively using
$\lambda (p)$ defined in \eqref{lmbda2}
and Kunita's inequality \eqref{eq:BDGPoisson2} for jump processes. 
Then in Theorems~\ref{t4.1} and \ref{t4.2} we respectively deal with
disease extinction and persistence in the mean for the system
\eqref{lsir1}-\eqref{lsir3}. 
We show that the threshold behavior of
the stochastic SIR system \eqref{lsir1}-\eqref{lsir3}
is determined by the basic reproduction number
\begin{equation}
  \label{rn}
\mathcal{\widebar{R}}_0=\mathcal{R}_0-\frac{\beta_2}{\mu+\varepsilon +\eta}
=
\frac{\beta\Lambda / \mu - \beta_2}{\mu+\varepsilon +\eta}
\end{equation}
 which differs from \cite{wulia} due to the quantity
$$
 \beta_2 :=\frac{1}{2}\varrho_{2,2} +\int_{\real^3 \setminus\{0\}}
  \left( \gamma_2(z)-\log (1+\gamma_2(z)) \right)
  \nu(dz).
$$
In Section~\ref{sec4} we present numerical simulations
based on tempered stable processes with parameter $\alpha \in (0,1)$.
We show in particular that the addition of a jump component to the system
\eqref{lsir1}-\eqref{lsir3} may result into the extinction of the
infected and recovered populations as $\alpha \in (0,1)$ becomes large enough
and the variance of random fluctuations increases, 
which is consistent with related observations in the literature,
see e.g. \cite{ycai}.
In addition, we note  
that this phenomenon can be observed when the noise variances are normalized to
identical values, showing that the shape of the distribution alone can affect the
long term behavior of the system. 
  The proofs
which are similar to the literature, see \cite{wulia},
are presented in the Appendix for completeness.
\section{Large time estimates}
\label{sec2}
 For $f$ an integrable function on $[0,t]$, we denote
  $$
  \langle f\rangle_t=\frac{1}{t}\int_0^tf(s)ds,
  \quad
  \langle f\rangle^\ast=\limsup\limits_{t\rightarrow\infty}\frac{1}{t}\int_0^tf(s)ds,
  \quad
  \langle f\rangle_\ast=\liminf\limits_{t\rightarrow\infty}\frac{1}{t}\int_0^tf(s)ds,
  \quad
 t>0.
  $$
 In addition, we assume that
 the jump coefficients $\gamma_i(z)$ in \eqref{sdjakl} 
 satisfy

 \bigskip

 \noindent
\textbf{$(H_1)$} ~$\displaystyle \int_{\real^3 \setminus\{0\}}|\gamma_i(z)|^2\nu(dz)<\infty$,
$i=1,2,3$,

\bigskip

\noindent
 together with the condition:
\\

\noindent
\textbf{$(H_2)$} ~$\gamma_i(z)>-1$, $\nu (dz)$-$a.e.$, and
$\displaystyle
\int_{\real^3 \setminus\{0\}}
\big(
\gamma_i(z) - \log ( 1 + \gamma_i(z))
\big)
\nu(dz)<\infty$, \ $i=1,2,3$.

  \medskip

\begin{prop}
  \label{Theorem 2.1} Under $(H_1)$-$(H_2)$,
  for any given initial data $(S_0,I_0,$ $R_0)\in\real _+^3$, 
  the system \eqref{lsir1}-\eqref{lsir3} admits a unique
  positive solution $(S_t,I_t,R_t)_{t\in \real_+}$
 which exists in $(0,\infty )^3$ for all $t\geq 0$, almost surely.
\end{prop}
\begin{Proof}
  By Theorem~6.2.11 in \cite{applebk2} or Theorem~2.1 in \cite{Bao1},
 the system \eqref{lsir1}-\eqref{lsir3}
 admits a unique local solution $(S_t,I_t,R_t)_{t\in(0,\tau_e]}$
 up to the explosion time $\tau_e$
 for any initial data $(S_0,I_0,R_0)\in\real _+^3$,
 since it has affine coefficients by $(H_1)$.
 In addition, by $(H_2)$ we have $\gamma_i(z)>-1$, $\nu (dz)$-$a.e.$, $i=1,2,3$,
 hence the solution is positive.
 The remaining of the proof
 follows the lines of proof of
 Theorem~1 in \cite{amllevy} by noting that
 the condition $(H_2)$ page~868 therein can be replaced by
 $(H_2)$ above.
\end{Proof}
Next, given $\lambda (p)$ defined in \eqref{lmbda2}
we let
$$
\Vert \varrho \Vert_\infty : = \max_{i=1,2,3} \sum_{j=1}^3 |\varrho_{i,j}|,
$$
and we consider the following condition:
\\

\textbf{$(H_3^{(p)})$}
~$\displaystyle
\mu > \frac{p-1}{2}\Vert \varrho \Vert_\infty + \frac{\lambda (p)}{p}, \qquad p>1$. 
\\

\noindent
  The proofs of Lemmas~\ref{l3.1}-\ref{l3.4} below present several
significant differences from the arguments of \cite{zhaoamc}
 and \cite{wulia}. 
First, our arguments do not require the finiteness
of the L\'evy measure $\nu (dz)$,
and in the proof of our Lemma~\ref{l3.1}
we replace the Burkholder-Davis-Gundy inequality
\eqref{bdg0} for
continuous processes used in \cite{wulia} with the simpler
bound \eqref{bd} below, 
as \eqref{bdg0} is not valid for jump processes and the inequality used at the beginning
of the proof of Lemma~2.1 in \cite{wulia} may not hold
in general because the compensated Poisson process $\tilde{N}(t)$
can have a negative drift.
Second, the proof of our Lemma~\ref{l3.4}
uses Kunita's inequality \eqref{eq:BDGPoisson2}
for jump processes
instead of relying on the Burkholder-Davis-Gundy inequality
\eqref{bdg0} for continuous processes. 

\bigskip

\noindent
 In the sequel, we consider the condition

 \bigskip

  \textbf{$(H_4^{(p)})$} ~ $\displaystyle 
 \int_{\real^3 \setminus\{0\}} | (1+\overline{\gamma}(z))^p-1 | \nu(dz) < \infty,
 \qquad p >1$. 

\medskip

\medskip

\begin{lemma}\label{l3.1}
  Assume that $(H_1)$-$(H_2)$ and $(H_3^{(p)})$-$(H_4^{(p)})$ hold for some $p>1$,
  and let $(S_t,I_t,R_t)$ be the solution of the system \eqref{lsir1}-\eqref{lsir3} with
   initial condition $(S_0,I_0,R_0)\in\real _+^3$.
Then we have
\begin{equation}\notag
\lim\limits_{t\rightarrow\infty}\frac{S_t}{t}=0,\quad \lim\limits_{t\rightarrow\infty}\frac{I_t}{t}=0,\quad
\lim\limits_{t\rightarrow\infty}\frac{R_t}{t}=0, \quad  \mathbb{P}\mbox{-}a.s.
\\
\end{equation}
\end{lemma}
\begin{Proof}
  We let $U_t:=S_t+I_t+R_t$ and
  $$
  H_t(z):=\gamma_1(z)S_t+\gamma_2(z)I_t+\gamma_3(z)R_t,
  \qquad
  z\in \real^3, \quad t\in \real_+,
  $$
  with the inequality $H_t(z)\leq \overline{\gamma} (z)U_t$,
  $z\in \real^3$.
  Applying the It\^{o} formula with jumps (see e.g. Theorem 1.16 in \cite{yiteng})
  to the function $x\mapsto V(x)=(1+x)^p$, we obtain
  \begin{eqnarray}
    \nonumber
    \lefteqn{
      dV(U_t)  =  p(1+U_t)^{p-1}(\Lambda-\mu U_t-\varepsilon  I_t)dt
    }
    \\
    \nonumber
    & &
    +\frac{p(p-1)}{2}(1+U_t)^{p-2}
    \big(\varrho_{1,1}S_t^2+\varrho_{2,2}I_t^2+\varrho_{3,3}R_t^2
    + 2 \varrho_{1,2} S_tI_t
    + 2 \varrho_{1,3} S_tR_t
    + 2 \varrho_{2,3} I_tR_t
    \big)dt
    \\
    \nonumber
& & +p(1+U_t)^{p-1} ( S_tdB^\varrho_1(t)+I_tdB^\varrho_2(t)+R_tdB^\varrho_3(t) ) \\
    \nonumber
& & +\int_{\real^3 \setminus\{0\}} ( (1+ U_{t}+H_t(z) ) ^p~-(1+U_{t})^p-p(1+U_{t})^{p-1}H_t(z) ) \nu(dz)dt\\
    \nonumber
& & +\int_{\real^3 \setminus\{0\}} ( (1+ U_{t^-}+H_{t^-}(z) ) ^p-(1+U_{t^-})^p ) \tilde{N}(dt,dz)
\\
\nonumber
&=& LV(U_t)dt+p(1+U_t)^{p-1} ( S_tdB^\varrho_1(t)+I_tdB^\varrho_2(t)+R_tdB^\varrho_3(t))
    \\
    \label{lv0}
    & & +\int_{\real^3 \setminus\{0\}} \big( (1+ U_{t^-}+H_{t^-}(z) ) ^p-(1+U_{t^-})^p \big) \tilde{N}(dt,dz),
\end{eqnarray}
 where we let
\begin{eqnarray*}
  \lefteqn{
    \!
    L V(U_t) := p(1+U_t)^{p-1}(\Lambda-\mu U_t-\varepsilon I_t)
  }
  \\
& &
     +\frac{p(p-1)}{2}(1+U_t)^{p-2}
    \big(\varrho_{1,1}S_t^2+\varrho_{2,2}I_t^2+\varrho_{3,3}R_t^2
    + 2 \varrho_{1,2} S_tI_t
    + 2 \varrho_{1,3} S_tR_t
    + 2 \varrho_{2,3} I_tR_t
    \big)
\\
& & +\int_{\real^3 \setminus\{0\}} ( (1+ U_{t}+H_t(z) ) ^p~-(1+U_{t})^p-p(1+U_{t})^{p-1}H_t(z) ) \nu(dz).
\end{eqnarray*}
 We note that for all $z\in \real^3$ and $t\in\real _+$ there exists
 $\theta \in(0,1)$ such that
  \begin{eqnarray*}
   \lefteqn{
 \! \! \! \! \! \! \! \!     (1+U_t+H_t(z))^p-(1+U_t)^p-p(1+U_t)^{p-1}H_t(z)
   }
   \\
   &=&(1+U_{t})^p+p(1+U_{t})^{p-1}H_t(z)
   +\frac{p(p-1)}{2} (1+ U_{t}+\theta H_t(z) ) ^{p-2}H_t^2(z)\\
   &&
    -(1+U_{t})^p-p(1+U_{t})^{p-1}H_t(z)\\
&=&\frac{p(p-1)}{2} (1+ U_{t}+\theta H_t(z)) ^{p-2}H_t^2(z)\\
    &\leq&\frac{p(p-1)}{2}\max ( 2^{p-3},1 ) ((1+ U_{t})^{p-2}+\theta H_t^{p-2}(z) ) H_t^2(z)
    \\
    &\leq&c_p(1+U_{t})^{p-2}H_t^2(z)+c_pH_t^p(z)\\
       &\leq&c_p (1+U_{t})^{p-2}U_{t}^2 ( \overline{\gamma}^2(z)+\overline{\gamma}^p(z) ),
       \quad
       z\in \real^3, \quad t\in\real _+,
\end{eqnarray*}
  where we used the bound $H_t(z)\leq \overline{\gamma}(z)U_t$,
  with $c_p:= p (p-1)\max ( 2^{p-3},1 )/2$.
  It then follows that
\begin{eqnarray}
\nonumber
\lefteqn{
   LV(U_t) \leq p(1+U_t)^{p-2}((1+U_t)(\Lambda-\mu U_t )-\varepsilon (1+U_t)I_t)
}
\\
   \nonumber
   & &
   +\frac{p(p-1)}{2}(1+U_t)^{p-2}
    \big(\varrho_{1,1}S_t^2+\varrho_{2,2}I_t^2+\varrho_{3,3}R_t^2
    + 2 \varrho_{1,2} S_tI_t
    + 2 \varrho_{1,3} S_tR_t
    + 2 \varrho_{2,3} I_tR_t
    \big)
   \\
   \nonumber
   & & +c_p (1+U_{t})^{p-2}U_{t}^2 \int_{\real^3 \setminus\{0\}} \overline{\gamma}^2(z) \nu(dz)
   +c_p (1+U_{t})^{p-2}U_{t}^2\int_{\real^3 \setminus\{0\}} \overline{\gamma}^p(z) \nu(dz)
   \\
   \nonumber
   &\leq & p(1+U_t)^{p-2}(-\mu U_t^2+(\Lambda-\mu)U_t+\Lambda)+
   \frac{p (p-1)}{2}(1+U_t)^{p-2} \Vert \varrho \Vert_\infty U_{t}^2
   \\
   \nonumber
   &&+\frac{pc_p}{p} (1+U_t)^{p-2}U_{t}^2\int_{\real^3 \setminus\{0\}} ( \overline{\gamma}^2(z)+\overline{\gamma}^p(z) ) \nu(dz)
   \\
   \label{lv}
&=& p(1+U_t)^{p-2}(- bU_t^2+(\Lambda-\mu)U_t+\Lambda),
\end{eqnarray}
 where
 $$
 b:=\mu-\frac{p-1}{2}\Vert \varrho \Vert_\infty - \frac{\lambda (p)}{p}>0
 $$
 by $(H_3^{(p)})$.
 Next, for any $k\in\real $ it holds that
 \begin{eqnarray}
  \nonumber 
  e^{kt}(1+U_t)^p&=& (1+U_0)^p +
  \int_0^t e^{ks} ( k(1+U_s)^p+LV(U_s) ) ds
  \\
  \nonumber
   & & + p \int_0^t e^{ks}(1+U_s)^{p-1}( S_sdB^\varrho_1(s)+I_sdB^\varrho_2(s)+R_sdB^\varrho_3(s))
   \\
   \nonumber
   & & +\int_0^t e^{ks}\int_{\real^3 \setminus\{0\}}\big( ( 1+U_{s^-}+H_{s^-}(z))^p-(1+U_{s^-})^p
   \big)
   \tilde{N}(ds,dz),
\end{eqnarray}
 hence by taking expectations in
 \eqref{lv0} and in view of \eqref{lv}, for any $k<bp$ we obtain
 \begin{eqnarray}
\nonumber 
   \lefteqn{
     e^{kt}\E[(1+U_t)^p]
    = 
   (1+U_0)^p+\E\left[
    \int_0^t e^{ks} \big( k(1+U_s)^p+LV(U_s) \big) ds\right].
   }
   \\
   \nonumber
   & \leq &
   (1+U_0)^p+\E\left[
     \int_0^t e^{ks} \left( k (1+U_s)^p+p (1+U_s)^{p-2} (\Lambda +(\Lambda-\mu)U_s - bU_s^2 )\right) ds\right]
   \\
   \nonumber
   & = &
   (1+U_0)^p+p \E\left[
     \int_0^t
     e^{ks}U_s^{p-2}\left( -\left(b-\frac{k}{p}\right)U_s^2
     +\left( \Lambda-\mu+\frac{2k}{p}\right)U_s+\Lambda+\frac{k}{p}
  \right)ds\right]
   \\
   \nonumber
   & \leq &
   (1+U_0)^p+pM   \int_0^t
     e^{ks}
     ds
  \\
  \nonumber
  & = &  (1+U_0)^p+\frac{pM}{k}e^{kt}, \quad t\in \real_+,
\end{eqnarray}
where
$$
0<M:= 1 + \sup\limits_{x\in\real _+}
x^{p-2}\left( -\left( b-\frac{k}{p}\right)x^2+
\left(
\Lambda-\mu+\frac{2k}{p}\right)x+\Lambda+\frac{k}{p}
\right)<\infty.
$$
 Hence, for any $k\in(0,bp)$ we have
\begin{equation}\notag
\limsup\limits_{t\rightarrow\infty}\E[(1+U_t)^p]\leq\frac{pM}{k},
\end{equation}
which implies that there exists $M_0>0$ such that
\begin{equation}\label{dvp6}
  \E[(1+U_t)^p]\leq M_0,\qquad t\in \real_+.
\end{equation}
 Now, by \eqref{lv0} and \eqref{lv} there will be
\begin{eqnarray*}
  \lefteqn{
    (1+U_t)^p-(1+U_{k\delta})^p 
    \leq p \int_{k\delta}^t(1+U_s)^{p-2}( \Lambda +(\Lambda-\mu)U_s - bU_s^2 ) ds
  }
  \\
     & & + p \int_{k\delta}^t(1+U_s)^{p-1} ( S_sdB^\varrho_1(s)+I_sdB^\varrho_2(s)+R_sdB^\varrho_3(s) )
  \\
& & +\int_{k\delta}^t\int_{\real^3 \setminus\{0\}} \big(
    ( 1+U_{s^-}+H_{s^-}(z))^p-(1+U_{s^-})^p \big)
    \tilde{N}(ds,dz), \qquad t\geq k\delta ,
\end{eqnarray*}
 from which it follows that
 \begin{eqnarray*}
     \sup\limits_{k\delta\leq t\leq(k+1)\delta}(1+U_t)^p
     & \leq & (1+U_{k\delta})^p
     + p \sup\limits_{k\delta\leq t\leq(k+1)\delta}\left|\int_{k\delta}^t(1+U_s)^{p-2}(
     \Lambda +(\Lambda-\mu)U_s - bU_s^2)ds\right|
   \\
   &&+ p \sup\limits_{k\delta\leq t\leq(k+1)\delta}\left|\int_{k\delta}^t(1+U_s)^{p-1} ( S_sdB^\varrho_1(s)+I_sdB^\varrho_2(s)+R_sdB^\varrho_3(s) )
   \right|\\
&&+\sup\limits_{k\delta\leq t\leq(k+1)\delta}\left|\int_{k\delta}^t\int_{\real^3 \setminus\{0\}}
\big( ( 1+U_{s^-}+H_{s^-}(z) )^p-(1+U_{s^-})^p \big) \tilde{N}(ds,dz)\right|.
\end{eqnarray*}
Taking expectation on both sides, we obtain
$$
  \E\bigg[\sup\limits_{k\delta\leq t\leq(k+1)\delta}(1+U_t)^p\bigg] \leq
  \E[(1+U_{k\delta})^p]+J_1+J_2+J_3
 \leq M_0+J_1+J_2+J_3,
$$
where, for some $c_3>0$,
\begin{eqnarray*}
J_1&:=& p \E\bigg[
  \sup\limits_{k\delta\leq t\leq(k+1)\delta}\left|\int_{k\delta}^t(1+U_s)^{p-2}
  ( \Lambda + (\Lambda-\mu)U_s - bU_s^2 )ds\right|\bigg]
    \\
    &\leq &
    c_3\E\bigg[\sup\limits_{k\delta\leq t\leq(k+1)\delta}
      \int_{k\delta}^t(1+U_s)^p ds \bigg]\\
&= & c_3\E\bigg[\int_{k\delta}^{(k+1)\delta}(1+U_s)^pds\bigg]\\
& \leq & c_3\delta \E\bigg[\sup\limits_{k\delta\leq t\leq(k+1)\delta}(1+U_t)^p\bigg],
\end{eqnarray*}
and, for some $c_4>0$,
\begin{eqnarray*}
J_2&:=& p \E\bigg[
  \sup\limits_{k\delta\leq t\leq(k+1)\delta}\left|\int_{k\delta}^t(1+U_s)^{p-1}
  \big(
  S_sdB^\varrho_1(s)+I_sdB^\varrho_2(s)+R_sdB^\varrho_3(s)
\big) \right|\bigg]
    \\
    &\leq & p \sqrt{32} 
    \E\bigg[
      \bigg(
        \int_{k\delta}^{(k+1)\delta}(1+U_s)^{2(p-1)}
      \big(
      \varrho_{1,1}S_s^2+\varrho_{2,2}I_s^2+\varrho_{3,3}R_s^2
      \big)
      ds
\bigg)^{1/2}
      \bigg]
  \\
    &\leq & p \sqrt{32 \delta \Vert \varrho \Vert_\infty}
    \E\bigg[
      \bigg(
      \sup\limits_{k\delta\leq t\leq(k+1)\delta}(1+U_t)^{2p}
      \bigg)^{1/2}
\bigg]
    \\
    & = &c_4 \sqrt{\delta} \E\bigg[\sup\limits_{k\delta\leq t\leq(k+1)\delta}(1+U_t)^p\bigg],
    \end{eqnarray*}
where we used the Burkholder-Davis-Gundy inequality \eqref{bdg0} for
continuous martingales, see Theorem~IV.4.48 page 193 
of \cite{protterb2005}, or Theorem~7.3
in Chapter~1 of \cite{mao2008}.
Furthermore, we have
\begin{eqnarray}
  \nonumber
  J_3&=&
  \E\left[
      \sup\limits_{k\delta\leq t\leq(k+1)\delta}\left|\int_{k\delta}^t\int_{\real^3 \setminus\{0\}}
      \big( (1+ U_{s^-}+H_{s^-}(z) ) ^p-(1+U_{s^-})^p\big) \tilde{N}(ds,dz)\right|\right]
  \\
 \nonumber
   & \leq &
  \E\left[
    \left|\int_{k\delta}^{(k+1)\delta}
    \int_{\real^3 \setminus\{0\}}
      \big( (1+ U_{s^-}+H_{s^-}(z) ) ^p-(1+U_{s^-})^p\big) {N}(ds,dz)\right|\right]
  \\
\nonumber
    & &
+   
  \E\left[
    \left|\int_{k\delta}^{(k+1)\delta}
    \int_{\real^3 \setminus\{0\}}
      \big( (1+ U_{s^-}+H_{s^-}(z) ) ^p-(1+U_{s^-})^p\big) ds \nu (dz)\right|\right]
\\
\nonumber
  & = &
2 \E\left[
    \left|\int_{k\delta}^{(k+1)\delta}
    \int_{\real^3 \setminus\{0\}}
      \big( (1+ U_{s^-}+H_{s^-}(z) ) ^p-(1+U_{s^-})^p\big) ds \nu (dz)\right|\right]
  \\
\nonumber
    &\leq &
2 
  \E\left[
    \int_{k\delta}^{(k+1)\delta}\int_{\real^3 \setminus\{0\}}
    (1+U_{s^-})^p
    | ( 1 + \overline{\gamma}(z)
     )^p-1
      |
      ds \nu (dz)
      \right]
  \\
\nonumber
      &= &
2 
    \E \left[ \int_{k\delta}^{(k+1)\delta}
    (1+U_s)^p ds \right]
    \int_{\real^3 \setminus\{0\}}
    | (1+\overline{\gamma}(z))^p-1 | \nu(dz)
    \\
    \label{bd}
    &\leq& 2 \delta 
    \E\Bigg[
      \sup\limits_{k\delta\leq t \leq(k+1)\delta}(1+U_t)^p \Bigg]
    \int_{\real^3 \setminus\{0\}}
    | (1+\overline{\gamma}(z))^p-1 | \nu(dz)
. 
\end{eqnarray}
 Therefore, we have
 \begin{eqnarray}
\label{Esup}
   \lefteqn{
     \E\bigg[ \sup\limits_{k\delta\leq t\leq(k+1)\delta}(1+U_t)^p\bigg]
   }
   \\
   \nonumber
   & \leq & \E[(1+U_{k\delta})^p]
     +\bigg( c_3\delta+c_4 \sqrt{\delta}+
     2 \delta 
     \int_{\real^3 \setminus\{0\}}
     | (1+\overline{\gamma}(z))^p-1 | \nu(dz)
          \bigg)
  \E\bigg[\sup\limits_{k\delta\leq t\leq(k+1)\delta}(1+U_t)^p\bigg].
   \end{eqnarray}
 Furthermore, from
 $(H_4^{(p)})$ we can choose $\delta>0$ such that
$$
c_3\delta+c_4 \sqrt{\delta}+
     2 \delta \int_{\real^3 \setminus\{0\}}
     | (1+\overline{\gamma}(z))^p-1 | \nu(dz)
<\frac{1}{2},
$$
 and, combining \eqref{dvp6} with \eqref{Esup}, we obtain
\begin{equation}\label{bofE}
  \E\bigg[\sup\limits_{k\delta\leq t\leq(k+1)\delta}(1+U_t)^p\bigg]\leq2\E[(1+U_{k\delta})^p]
  \leq 2M_0.
\end{equation}
Let now $\varepsilon >0$ be arbitrary.
By Chebyshev's inequality, we get
\begin{equation}\notag
  \mathbb{P}\bigg(
  \sup\limits_{k\delta\leq t\leq(k+1)\delta}(1+U_t)^p>(k\delta)^{1+\varepsilon }
  \bigg)
  \leq\frac{1}{(k\delta)^{1+\varepsilon }}
  \E\bigg[\sup\limits_{k\delta\leq t\leq(k+1)\delta}(1+U_t)^p\bigg]
  \leq\frac{2M_0}{(k\delta)^{1+\varepsilon }}
\end{equation}
for all $k \geq 1$.
Then, by the Borel-Cantelli lemma (see Lemma 2.4 in Chapter~1
of \cite{mao2008}) it follows that for almost all $\omega\in\Omega$, the bound
\begin{equation}
  \nonumber 
  \sup\limits_{k\delta\leq t\leq(k+1)\delta}(1+U_t)^p\leq(k\delta)^{1+\varepsilon },
\end{equation}
holds for all but finitely many $k$. Thus,
for almost all $\omega\in\Omega$
there exists $k_0(\omega)$ such that whenever $k\geq k_0 (\omega)$
we have
$$
\frac{\log (1+U_t)^p}{\log t}
\leq
1+\varepsilon, \qquad \varepsilon >0,
\quad k\delta\leq t\leq(k+1)\delta,
$$
 hence
\begin{equation}\notag
  \limsup\limits_{t\rightarrow\infty}\frac{\log U_t}{\log t}\leq\limsup\limits_{t\rightarrow\infty}\frac{\log (1+U_t)}{\log t}\leq\frac{1}{p}, \quad  \mathbb{P}\mbox{-}a.s., \quad
  p>1.
\end{equation}
In other words,
for any $\xi \in (0,1-1/p)$
there exists an a.s. finite random time $\widebar{T} (\omega)$
such that 
\begin{equation}\notag
  \log U_t\leq \left( \frac{1}{p}+\xi \right)\log t,
  \qquad
t\geq \widebar{T}.
\end{equation}
It follows that
\begin{equation}\notag
  \limsup\limits_{t\rightarrow\infty}\frac{U_t}{t}\leq\limsup\limits_{t\rightarrow\infty}\frac{t^{\xi + 1/p}}{t}=0,
\end{equation}
therefore we have
\begin{equation}\nonumber 
  \limsup\limits_{t\rightarrow\infty}\frac{S_t}{t}\leq0,\quad \limsup\limits_{t\rightarrow\infty}\frac{I_t}{t}\leq0,
\quad \limsup\limits_{t\rightarrow\infty}\frac{R_t}{t}\leq0, \quad  \mathbb{P}\mbox{-}a.s.
\end{equation}
This, together with the positivity of the solution, implies that
\begin{equation}\notag
\lim\limits_{t\rightarrow\infty}\frac{S_t}{t}=0,\quad \lim\limits_{t\rightarrow\infty}\frac{I_t}{t}=0,\quad \lim\limits_{t\rightarrow\infty}\frac{R_t}{t}=0, \quad  \mathbb{P}\mbox{-}a.s.
\end{equation}
\end{Proof}
The next Lemma~\ref{l3.4} is proved by using
Kunita's inequality \eqref{eq:BDGPoisson2}
for jump processes instead of
the Burkholder-Davis-Gundy inequality
\eqref{bdg0} for continuous martingales. 
\begin{lemma}
  \label{l3.4}
    Assume that $(H_1)$-$(H_2)$ and $(H_3^{(p)})$-$(H_4^{(p)})$ hold for some $p>2$,
  and let $(S_t,I_t,R_t)$ be the solution of the system \eqref{lsir1}-\eqref{lsir3} with
  initial condition $(S_0,I_0,R_0)\in\real _+^3$.
  Then we have
$$
  \lim\limits_{t\rightarrow\infty}\frac{1}{t}\int_0^t S_{r^-}\int_{\real^3 \setminus\{0\}}\gamma_1(z)\tilde{N}(dr,dz)=0,\quad
  \lim\limits_{t\rightarrow\infty}\frac{1}{t}\int_0^t I_{r^-}\int_{\real^3 \setminus\{0\}}\gamma_2(z)\tilde{N}(dr,dz)=0,
  $$ and
\begin{equation}\notag
\lim\limits_{t\rightarrow\infty}\frac{1}{t}\int_0^t R_{r^-} \int_{\real^3 \setminus\{0\}}\gamma_3(z)\tilde{N}(dr,dz)=0, \quad  \mathbb{P}\mbox{-}a.s.\end{equation}
\end{lemma}
\begin{Proof}
  Denote
  $$
  X_1(t):=\int_0^tS_{r^-}\int_{\real^3 \setminus\{0\}}\gamma_1(z)\tilde{N}(dr,dz),
  \quad
  X_2(t):=\int_0^tI_{r^-}\int_{\real^3 \setminus\{0\}}\gamma_2(z)\tilde{N}(dr,dz),
  $$
  and
  $$
  X_3(t):=\int_0^tR_{r^-}\int_{\real^3 \setminus\{0\}}\gamma_3(z)\tilde{N}(dr,dz),
  \quad
  t\in \real_+.
  $$
  By Kunita's inequality \eqref{eq:BDGPoisson2}, for any $p \geq 2$
  there exists a positive constant $C_p$ such that
\begin{align*}
 &     \E\bigg[\sup\limits_{0< r\leq t}|X_1(r)|^p\bigg]
  \\
  &   \leq
  C_p\E\bigg[\bigg(\int_0^t|S_r|^2 \int_{\real^3 \setminus\{0\}}|\gamma_1(z)|^2\nu(dz)dr\bigg)^{p/2}\bigg]
  +C_p\E\bigg[\int_0^t |S_r|^p \int_{\real^3 \setminus\{0\}}|\gamma_1(z)|^p\nu(dz)dr\bigg]
  \\
  &=  C_p
  \bigg( \int_{\real^3 \setminus\{0\}}\gamma_1^2(z)\nu(dz) \bigg)^{p/2}
  \E \bigg[\bigg(\int_0^t|S_r|^2dr\bigg)^{p/2}\bigg]
   +C_p \left(
  \int_{\real^3 \setminus\{0\}}\gamma_1^p(z)\nu(dz) \right)
  \E\bigg[\int_0^t|S_r|^p dr\bigg]
  \\
  &\leq  C_p
  t^{p/2}
  \left( \int_{\real^3 \setminus\{0\}}\gamma_1^2(z)\nu(dz)\right)^{p/2}\E\bigg[\left(\sup\limits_{0< r\leq t}|S_r|^2\right)^{p/2}\bigg]
  \\
  & \quad
  +C_p \int_{\real^3 \setminus\{0\}}\gamma_1^p(z)\nu(dz)
  \int_0^t\E[|S_r|^p]dr
  \\
  &\leq  C_p
  t^{p/2}
  \left( \int_{\real^3 \setminus\{0\}}\gamma_1^2(z)\nu(dz) \right)^{p/2}
  \E\left[\sup\limits_{0< r\leq t}|S_r|^p\right]
   +C_pM_0t\int_{\real^3 \setminus\{0\}}\gamma_1^p(z)\nu(dz),
\end{align*}
where the bound \eqref{dvp6} has been used in the last inequality.
Combining the above inequality with \eqref{bofE} yields
\begin{eqnarray*}
  \lefteqn{
    \E\bigg[
      \sup\limits_{k\delta\leq t\leq(k+1)\delta}|X_1(t)|^p\bigg]
  }
  \\
   &\leq & C_p
  ( (k+1)\delta )^{p/2}
  \left(
  \int_{\real^3 \setminus\{0\}}\gamma_1^2(z)\nu(dz) \right)^{p/2}
  \E\bigg[\sup\limits_{k\delta\leq t\leq(k+1)\delta}|S_r|^p\bigg]
  \\
   & & +C_p M_0(k+1)\delta \int_{\real^3 \setminus\{0\}}\gamma_1^p(z)\nu(dz) \\
  &\leq & C_p2M_0 ( (k+1)\delta ) ^{p/2}
    \left( \int_{\real^3 \setminus\{0\}}\gamma_1^2(z)\nu(dz) \right)^{p/2}
 +C_p \delta M_0(k+1) \int_{\real^3 \setminus\{0\}}\gamma_1^p(z)\nu(dz).
\end{eqnarray*}
Let $\varepsilon >0$ be arbitrary. It follows from Doob's martingale inequality
(see e.g. Theorem~3.8 in Chapter~1 of \cite{mao2008}) that
\begin{eqnarray*}
  \lefteqn{
    \mathbb{P}\bigg( \sup\limits_{k\delta\leq t\leq(k+1)\delta}|X_1(t)|^p>(k\delta)^{1+\varepsilon +p/2}\bigg)
   \leq (k\delta)^{-1-\varepsilon -p/2}
  \E\bigg[ \sup\limits_{k\delta\leq t\leq(k+1)\delta}|X_1(t)|^p\bigg]
  }
  \\
    & \leq&\frac{C_p2M_0 ( (k+1)\delta )^{p/2}}{(k\delta)^{1+\varepsilon +p/2}}
      \left(\int_{\real^3 \setminus\{0\}}\gamma_1^2(z)\nu(dz) \right)^{p/2}
        +\frac{C_pM_0(k+1)\delta}{(k\delta)^{1+\varepsilon +p/2}}
      \int_{\real^3 \setminus\{0\}}\gamma_1^p(z)\nu(dz).
\end{eqnarray*}
By the Borel-Cantelli lemma it follows that for almost all $\omega\in\Omega$ the bound
\begin{equation}
  \nonumber 
  \sup\limits_{k\delta\leq t\leq(k+1)\delta}|X_1(t)|^p\leq(k\delta)^{1+\varepsilon +p/2}
\end{equation}
holds for all but finitely many $k$.
Thus, for almost all $\omega\in\Omega$
there exists $k_0 (\omega)$ such that for
all $k\geq k_0 (\omega)$ we have
\begin{equation}\notag
 \frac{\log |X_1(t)|}{\log t}
  \leq
  \frac{1}{2} + \frac{1+\varepsilon}{p},
  \qquad
  \varepsilon >0,
  \quad
  k\delta\leq t\leq(k+1)\delta,
\end{equation}
hence 
\begin{equation}\notag
\limsup_{t\to \infty} \frac{\log |X_1(t)|}{\log t}
  \leq
  \frac{1}{2} + \frac{1}{p}, \qquad p>2, 
\end{equation}
 and as in the proof of Lemma~\ref{l3.1}, this yields 
\begin{equation}\notag
  \limsup\limits_{t\rightarrow\infty}\frac{|X_1(t)|}{t}\leq\limsup\limits_{t\rightarrow\infty}\frac{t^{ 1/2+1/p}}{t}=0
  \end{equation}
since $p>2$, which shows that
\begin{equation}\notag
\lim\limits_{t\rightarrow\infty}\frac{|X_1(t)|}{t}=0, \quad  \mathbb{P}\mbox{-}a.s.\end{equation}
By similar arguments, we also obtain
\begin{equation}\notag
\lim\limits_{t\rightarrow\infty}\frac{X_2(t)}{t}=0,\quad \lim\limits_{t\rightarrow\infty}\frac{X_3(t)}{t}=0, \quad  \mathbb{P}\mbox{-}a.s.\end{equation}
\end{Proof}
\begin{remark}
  \label{r0} 
\noindent We note that
by the continuity of $p\mapsto {\lambda}(p)$
in \eqref{lmbda2},
in order 
for $(H_3^{(p)})$ to hold for some $p>2$ it suffices that
$(H_3^{(2)})$ be satisfied, i.e.
\vspace*{0.2cm}
\begin{equation}\nonumber 
\textbf{$(H_3^{(2)})$}: ~~\mu >\frac{\Vert \varrho \Vert_\infty }{2}+\frac{{\lambda}(2)}{2}.
\end{equation}
\end{remark}
The next Lemma~\ref{l3.3} can be proved
on \eqref{dvp6},
by noting that the argument of Lemma~2.2 in \cite{wulia} 
is valid for correlated Brownian motions
$(B^\varrho_1 (t), B^\varrho_2 (t), B^\varrho_3(t))$,
without requiring the continuity of
$(S_t,I_t,R_t)_{t\in \real_+}$.
\begin{lemma}
  \label{l3.3}
  Assume that $(H_1)$-$(H_2)$ and $(H_3^{(p)})$-$(H_4^{(p)})$ hold for some $p>1$,
  and let $(S_t,I_t,R_t)$ be the solution of \eqref{lsir1}-\eqref{lsir3} with
  initial condition
  $(S_0,I_0,R_0)\in\real _+^3$.
    Then, $\mathbb{P}$-$a.s$, we have
\begin{equation}\notag
\lim\limits_{t\rightarrow\infty}\frac{1}{t}\int_0^t S_rdB^\varrho_1(r)=0,~\lim\limits_{t\rightarrow\infty}\frac{1}{t}\int_0^t I_rdB^\varrho_2(r)=0,~
\lim\limits_{t\rightarrow\infty}\frac{1}{t}\int_0^t R_rdB^\varrho_3(r)=0.\end{equation}
\end{lemma}

\section{Extinction and persistence}
By virtue of the large time estimates
for the solution of \eqref{lsir1}-\eqref{lsir3} and its
diffusion and jump components
obtained in Section~\ref{sec2},
in this section
we determine the threshold behavior of the stochastic SIR epidemic model.
\\

In Theorems~\ref{t4.1} and \ref{t4.2} below,  
the extinction and persistence of the disease
 is characterized by means of the critical reproduction number
$\mathcal{\widebar{R}}_0$ in \eqref{rn},
  which shows that the
  additional environmental noise induced by L\'evy jumps
  can limit the outbreak of the disease.
    In the sequel, we let
  \begin{equation}
  \nonumber 
 0 < \beta_i :=\frac{1}{2}\varrho_{i,i}+\int_{\real^3 \setminus\{0\}}
  \left( \gamma_i(z)-\log (1+\gamma_i(z)) \right)
  \nu(dz),\quad i=1,2,3,
 \end{equation}
 which is finite under $(H_1)$, 
 and we consider the following condition:

 \vspace*{0.2cm}\textbf{$(H_5)$} ~$\displaystyle
\int_{\real^3 \setminus\{0\}} \big( \log (1+\gamma_i(z)) \big)^2\nu(dz)<\infty$.

\bigskip
 
\noindent
We note that the basic reproduction number
$  \mathcal{\widebar{R}}_0$
becomes lower in the presence of jumps.

\begin{theorem}\label{t4.1}
  {\em (Extinction)}.
  Assume that
  $(H_1)$-$(H_2)$, $(H_3^{(p)})$-$(H_4^{(p)})$ and $(H_5)$ hold for some $p>2$.
    If in addition
\begin{equation}\notag
  \mathcal{\widebar{R}}_0 :=\mathcal{R}_0-\frac{\beta_2}{\mu+\varepsilon +\eta} < 1,
  \end{equation}
 then for any initial condition
 $(S_0,I_0,R_0)\in\real _+^3$,
 the disease vanishes with probability one in large time,
 i.e. the solution $(S_t,I_t,R_t)$
 of \eqref{lsir1}-\eqref{lsir3} satisfies
\begin{equation}\notag
  \lim\limits_{t\rightarrow\infty}\langle S\rangle_t=\frac{\Lambda}{\mu},
\quad
\lim\limits_{t\rightarrow\infty}I_t=0
\mbox{ ~and~ }
\lim\limits_{t\rightarrow\infty}R_t=0, \quad  \mathbb{P}\mbox{-}a.s.
\end{equation}
\end{theorem}
  The proof of Theorem~\ref{t4.1}
  follows the lines of the proof of Theorem~2.1 in \cite{wulia},
  up to the new Condition~$(H_3^{(p)})$
  which allows for infinite L\'evy measures
  in Lemma~\ref{l3.1}. 
  For reference, the proof of Theorem~\ref{t4.1} is stated in
  Appendix.
\\

  Next, we consider the persistence of the system \eqref{lsir1}-\eqref{lsir3}.
 We recall that the system \eqref{lsir1}-\eqref{lsir3} is said to be persistent in the mean if
\begin{equation}\notag
\liminf\limits_{t\rightarrow\infty}\frac{1}{t}\int_0^tS_rdr>0,\quad \liminf\limits_{t\rightarrow\infty}\frac{1}{t}\int_0^tI_rdr>0,
\quad \liminf\limits_{t\rightarrow\infty}\frac{1}{t}\int_0^tR_rdr>0, \quad  \mathbb{P}\mbox{-}a.s.\end{equation}
In Theorem~\ref{t4.2}
we explore the conditions for the disease to be endemic, in other words,
sufficient conditions for the persistence of
the infected population $I_t$.
  \begin{theorem}\label{t4.2}
    {\em (Persistence)}.
    Assume that
    $(H_1)$-$(H_2)$, $(H_3^{(p)})$-$(H_4^{(p)})$ and $(H_5)$ hold
    for some $p>2$. 
    If in addition
    \begin{equation}\notag
  \mathcal{\widebar{R}}_0 :=\mathcal{R}_0-\frac{\beta_2}{\mu+\varepsilon+\eta} > 1,
  \end{equation}
 then for any initial condition
  $(S_0,I_0,R_0)\in\real _+^3$,
  the solution $(S_t,I_t,R_t)$ of \eqref{lsir1}-\eqref{lsir3}
  satisfies
  $$
  \lim\limits_{t\rightarrow\infty}\langle S\rangle_t=S^\ast+\frac{\beta_2}{\beta}, \quad
  \lim\limits_{t\rightarrow\infty}\langle I\rangle_t=\frac{\mu}{\beta}(\mathcal{\widebar{R}}_0-1),
  \quad
\lim\limits_{t\rightarrow\infty}\langle R\rangle_t=\frac{\eta}{\beta}(\mathcal{\widebar{R}}_0-1), \quad  \mathbb{P}\mbox{-}a.s.,
$$
where $S^\ast := ( \mu+\varepsilon+\eta ) / \beta$
is the equilibrium value for the
susceptible population $S_t$ in the corresponding deterministic SIR model. 
\end{theorem}
  For reference, the proof of Theorem~\ref{t4.2} is stated in
  Appendix.
 It follows the lines of the proof of Theorem~3.1 in \cite{wulia},
  up to the use of Lemma~\ref{Lemma 2.2} (in Appendix)
  which extends Lemma~2 of \cite{mazhien} to discontinuous functions.
\section{Numerical experiments} 
\label{sec4}
In this section, we provide numerical simulations for the behavior of 
\eqref{lsir1}-\eqref{lsir3}
using tempered stable processes.
The (compensated) one-dimensional tempered stable L\'evy process
$$
Y(t) = \int_0^t \int_{\real \setminus\{0\}} z \tilde{N}(ds,dz),
\qquad t\in \real_+,
$$
is defined by its L\'evy measure \eqref{measure2.0} on $\real \setminus \{0\}$ 
where $k_-,k_+,\lambda_-,\lambda_+ > 0 $ and $\alpha \in (0,2)$,
i.e.
\begin{equation}
\label{measure3} 
 \nu (dz )
 =
   \frac{k_-}{z^{\alpha +1}}
 \re^{-\lambda_- z } dz
 +
\frac{k_+}{z^{\alpha +1}} \re^{-\lambda_+ z } dz. 
\end{equation} 
As $\nu (\real )=\infty$, the tempered stable process
$(Y(t))_{t\in \real_+}$ is not covered by the proof arguments of \cite{amllevy},
\cite{wulia} and \cite{xiaobing}, in particular
 the quantity defined in
 \eqref{lmbda} is not finite in this case.

\subsubsection*{Random simulations}

We use the simulation algorithm of
\cite{rosinski4} for the tempered stable process with L\'evy measure \eqref{measure2}.
Consider $(\epsilon_j)_{j\geq 1}$ an
independent and identically distributed (i.i.d.) 
Bernoulli $( -\lambda_-, \lambda_+)$-valued 
random sequence with distribution $(k_-/(k_-+k_+), k_+/(k_-+k_+))$, 
$(\xi_j)_{j\geq 1}$ an i.i.d. uniform $U(0,1)$ random sequence, 
and $(\eta_j)_{j\geq 1}$, $(\eta^\prime_j)_{j\geq 1}$ 
 i.i.d. exponentially distributed random sequences with parameter $1$,
with $\Gamma_j:=\eta^\prime_1+\cdots+\eta^\prime_j$, $j\geq 1$.
We also let $(u_j)_{j\geq 1}$
denote an i.i.d. sequence of uniform random variables on $[0,T]$, where $T>0$,
and assume that the sequence $(\epsilon_j)_{j\geq 1}$,
$(\xi_j)_{j\geq 1}$,
$(\eta_j)_{j\geq 1}$,
$(\eta^\prime_j)_{j\geq 1}$,
and $(u_j)_{j\geq 1}$ are mutually independent. 
By Theorem~5.3 in \cite{rosinski4}, the 
 tempered stable process $Y (t)$
 with L\'evy measure \eqref{measure2}
 admits the following representations.
\\
\noindent
$(i)$
If $\alpha\in(0,1)$, set
\begin{equation}
  \nonumber 
  Y(t)=\sum_{j=1}^\infty\mathbb{I}_{(0,t]}(u_j)\min\left\{\left(\frac{T(k_-+k_+)}{\alpha\Gamma_j}\right)^{1/\alpha}
 \hskip-0.2cm
 ,
 \frac{\eta_j}{|\epsilon_j|} \xi_j^{1/\alpha} 
\right\}\frac{\epsilon_j}{|\epsilon_j|},
\qquad t\in[0,T].
\end{equation}
\noindent
$(ii)$
 If $\alpha\in(1,2)$, set
\begin{equation}
   \nonumber 
   Y(t)
 =  \sum_{j=1}^\infty\left(\mathbb{I}_{(0,t]}(u_j)
  \min\left\{\left(\frac{k_-+k_+}{\alpha\Gamma_j/T}\right)^{1/\alpha}
  \hskip-0.2cm
  ,
\frac{\eta_j }{|\epsilon_j|} \xi_j^{1/\alpha} 
\right\}\frac{\epsilon_j}{|\epsilon_j|}-
x_0 \frac{t}{T}\left(\frac{k_-+k_+}{\alpha j/T}\right)^{1/\alpha}\right)+tb_T,
\end{equation} 
 $t\in[0,T]$, with 
$x_0= ( k_--k_+ ) / (k_-+k_+)$, $x_1=k_+ \lambda_+^{-1-\alpha}-k_-\lambda_-^{-1-\alpha}$, and 
$$
b_T:= \frac{x_0}{T} \zeta\left(
\frac{1}{\alpha} \right)
\left(
\frac{T(k_-+k_+)}{\alpha}
\right)^{1/\alpha}
- x_1 \Gamma(1-\alpha),
$$
where $\zeta$ is the Riemann zeta function.

\begin{figure}[H]
\centering
\hskip-0.2cm
\begin{subfigure}{.48\textwidth}
\centering
\includegraphics[width=1.\textwidth]{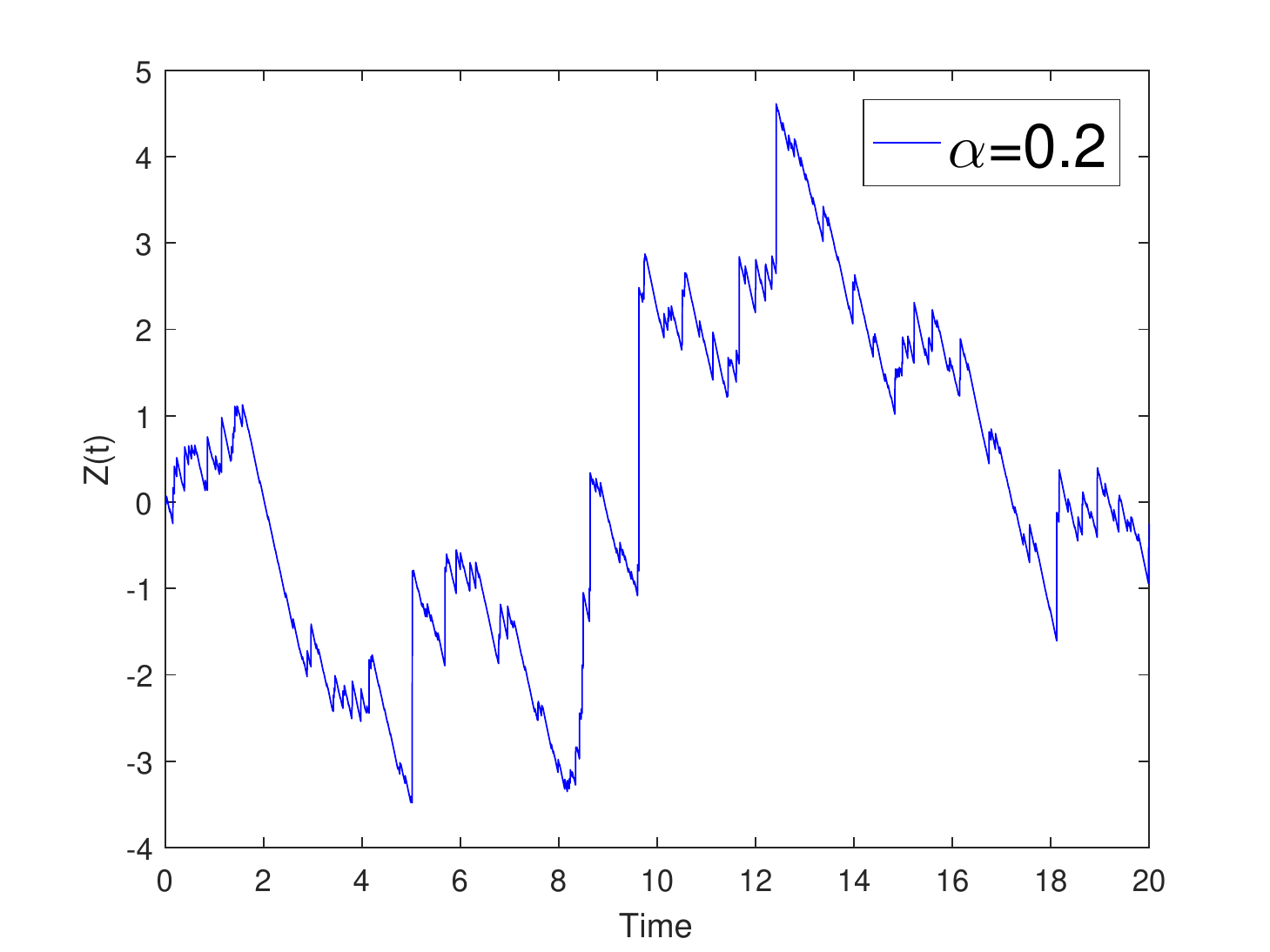}
\caption{\small Tempered stable process with $\alpha =0.2$} 
\end{subfigure}
\hskip0.5cm
\begin{subfigure}{.48\textwidth}
\centering
\includegraphics[width=1.\textwidth]{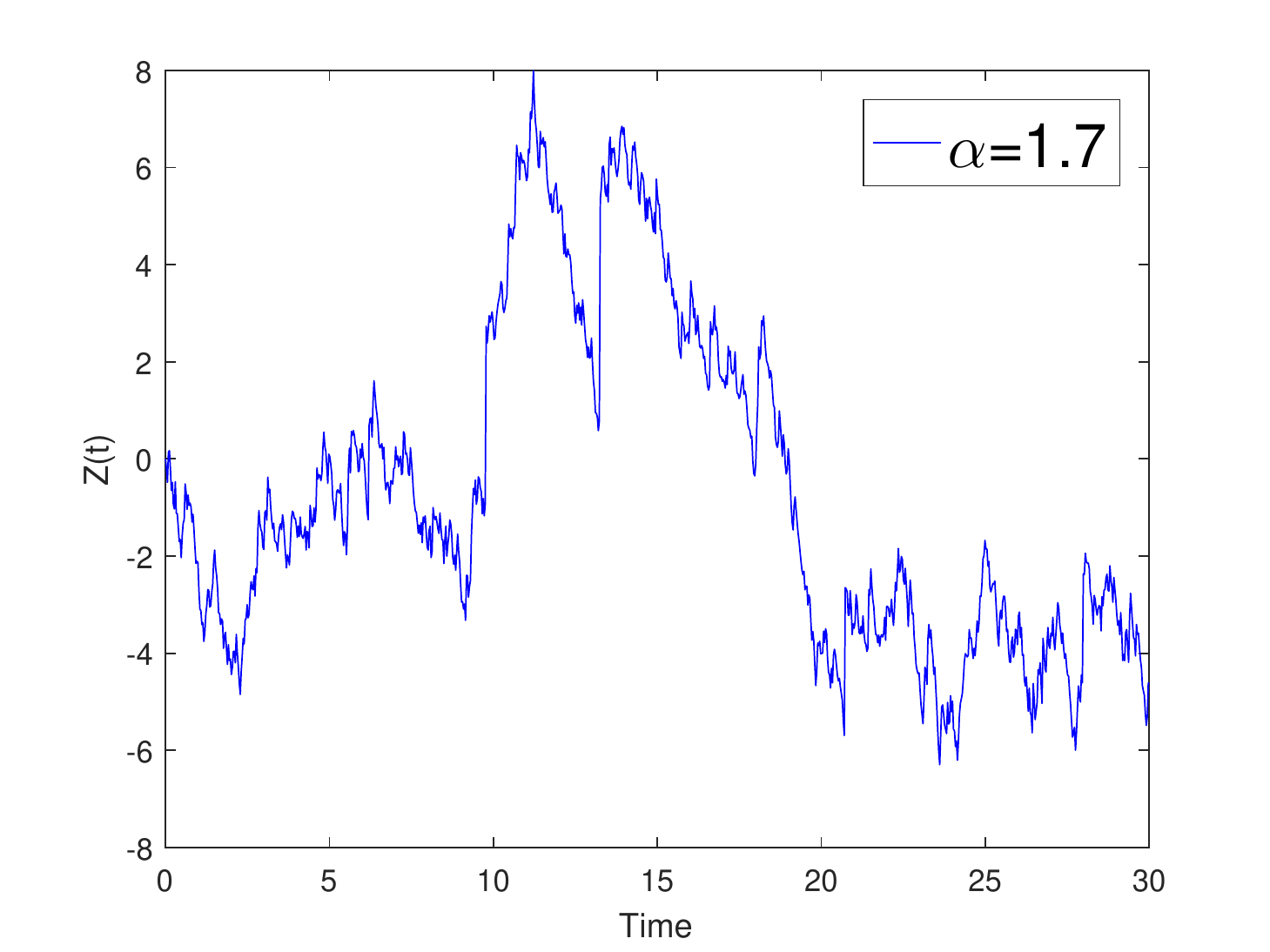}
\caption{\small Tempered stable process with $\alpha =1.7$} 
\end{subfigure}
\caption{\small Simulated sample paths of the tempered stable process.} 
\end{figure}

\noindent
 Next, we take
 $\gamma_i(z):= 
 \sigma_i z_i$
 with $\sigma_i>0$, $i=1,2,3$, and we consider the system
\begin{subequations}
\begin{empheq}[left=\empheqlbrace]{align}
\nonumber 
dS_t&= ( \Lambda-\beta S_tI_t-\mu S_t ) dt+S_tdB^\varrho_1 (t)+ \sigma_1 S_tdY(t),
\\[7pt]
\nonumber 
dI_t&= ( \beta S_tI_t-(\mu+\varepsilon+\eta)I_t ) dt+I_tdB^\varrho_2 (t)+ \sigma_2 I_tdY(t),
\\[7pt]
\nonumber 
    dR_t&= ( \eta I_t-\mu R_t ) dt+R_tdB^\varrho_3 (t)+ \sigma_3 R_tdY(t),
\end{empheq}
\end{subequations}
 i.e. \eqref{sdjakl} reads $Z_i(t) = B^\varrho_i (t)+Y(t)$, $i =1,2,3$, 
 and $(Y(t))_{t\in \real_+}$ is a one-sided tempered stable process 
 with $k_-=0$ in \eqref{measure2}. 
 We note that $(H_1)$-$(H_2)$, $(H_5)$ are
 satisfied, and that $(H_4^{(p)})$ holds for all $\alpha \in (0,1)$ and $p>1$. 
 In addition, letting $\overline{\sigma} := \max ( \sigma_1,\sigma_2,\sigma_3)$,
  the quantity
\begin{eqnarray*}
  \lambda (p) &  = &
  c_p \overline{\sigma}^2 \int_{\real^3 \setminus\{0\}}
  z^2 {\nu}(dz)
  +
  c_p \overline{\sigma}^p \int_{\real^3 \setminus\{0\}}
      z^p {\nu}(dz)
\\
& = & c_p k_+ \overline{\sigma}^2 
\frac{\Gamma(2-\alpha)}{\lambda_+^{2-\alpha}}
+
c_p k_+ \overline{\sigma}^2
\frac{\Gamma(p-\alpha)}{\lambda_+^{p-\alpha}}
\end{eqnarray*} 
 in \eqref{lmbda2} is finite when $p > \alpha$,
 where $\Gamma(\cdot)$ is the Gamma function.
 We note that the variance 
 $k_+ t \Gamma(2-\alpha)/\lambda_+^{2-\alpha}$
 of the one-sided tempered stable process
 $Y(t)$ is an increasing function of $\alpha \in (0,1)$
 when $\lambda_+\geq 1$.

 \bigskip
 
 \noindent
 First, we take $\alpha=0.7$, $\lambda_+=1.2$, $k_+=2.8$ with the
initial condition $(S_0,I_0,R_0)=(1.6,0.4,0.04)$
and the parameters
$\Lambda=8$, $\mu=5.3$, $\beta=4.8$, $\eta=1$ and $\varepsilon=0.5$. 
The covariance matrix is set at
$\varrho=10^{-2} \left(\begin{array}{ccc}
  4  & 3.2   & 3.0 \\
  3.2 & 4    & 3.84 \\
  3   & 3.84 & 4.69
\end{array}\right)$,
with $\sigma_1=0.2$, $\sigma_2=0.8$ and $\sigma_3=0.5$,
in which case Condition~$(H_3^{(2)})$ is also satisfied by 
Remark~\ref{r0}. 

\vskip-0.2cm

\begin{figure}[H]
\centering
\hskip-0.2cm
\begin{subfigure}{.5\textwidth}
\centering
\includegraphics[width=1\textwidth]{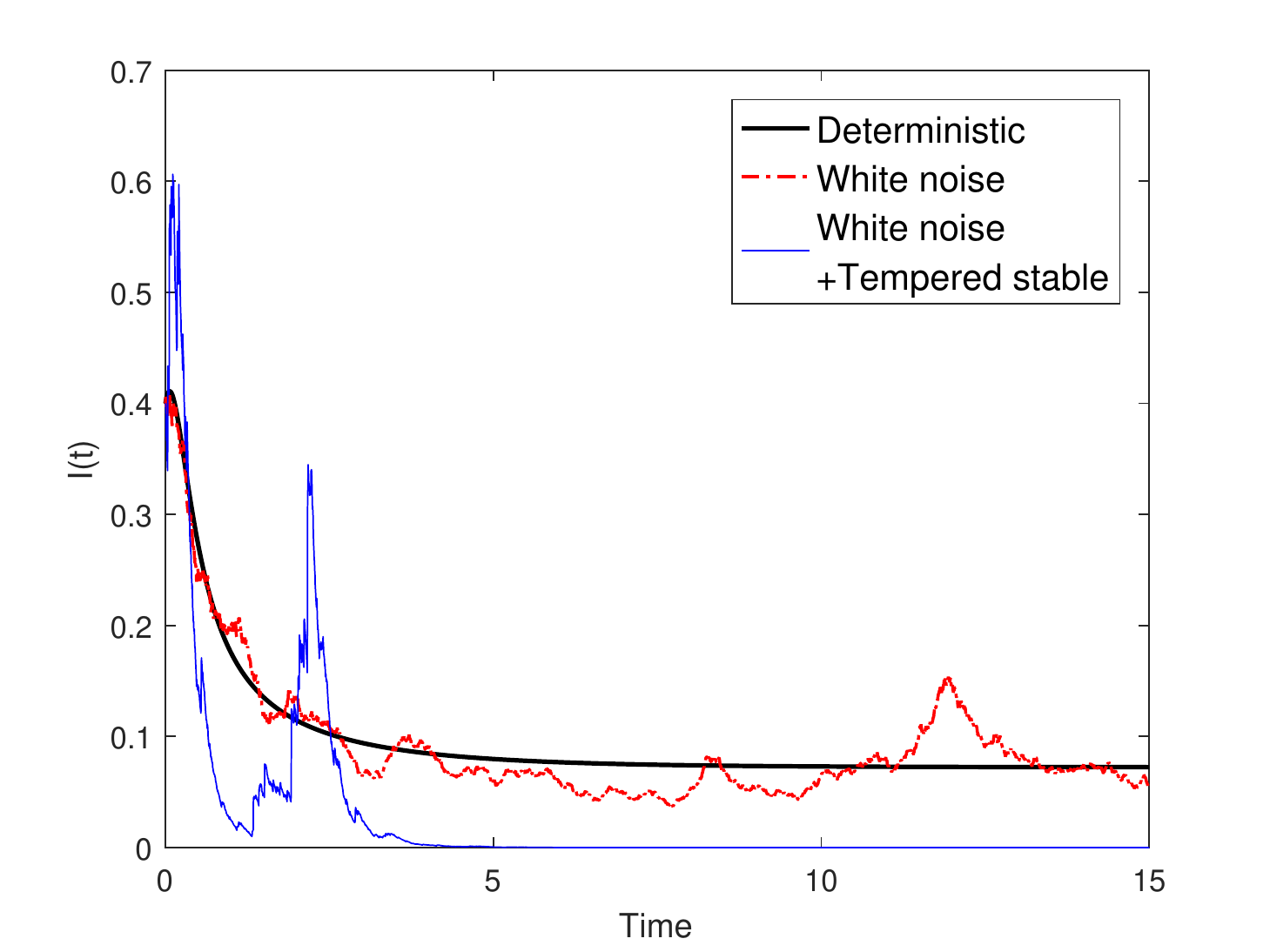}
\caption{\small Extinction of the infected population} 
\end{subfigure}
\begin{subfigure}{.5\textwidth}
\centering
\includegraphics[width=1\textwidth]{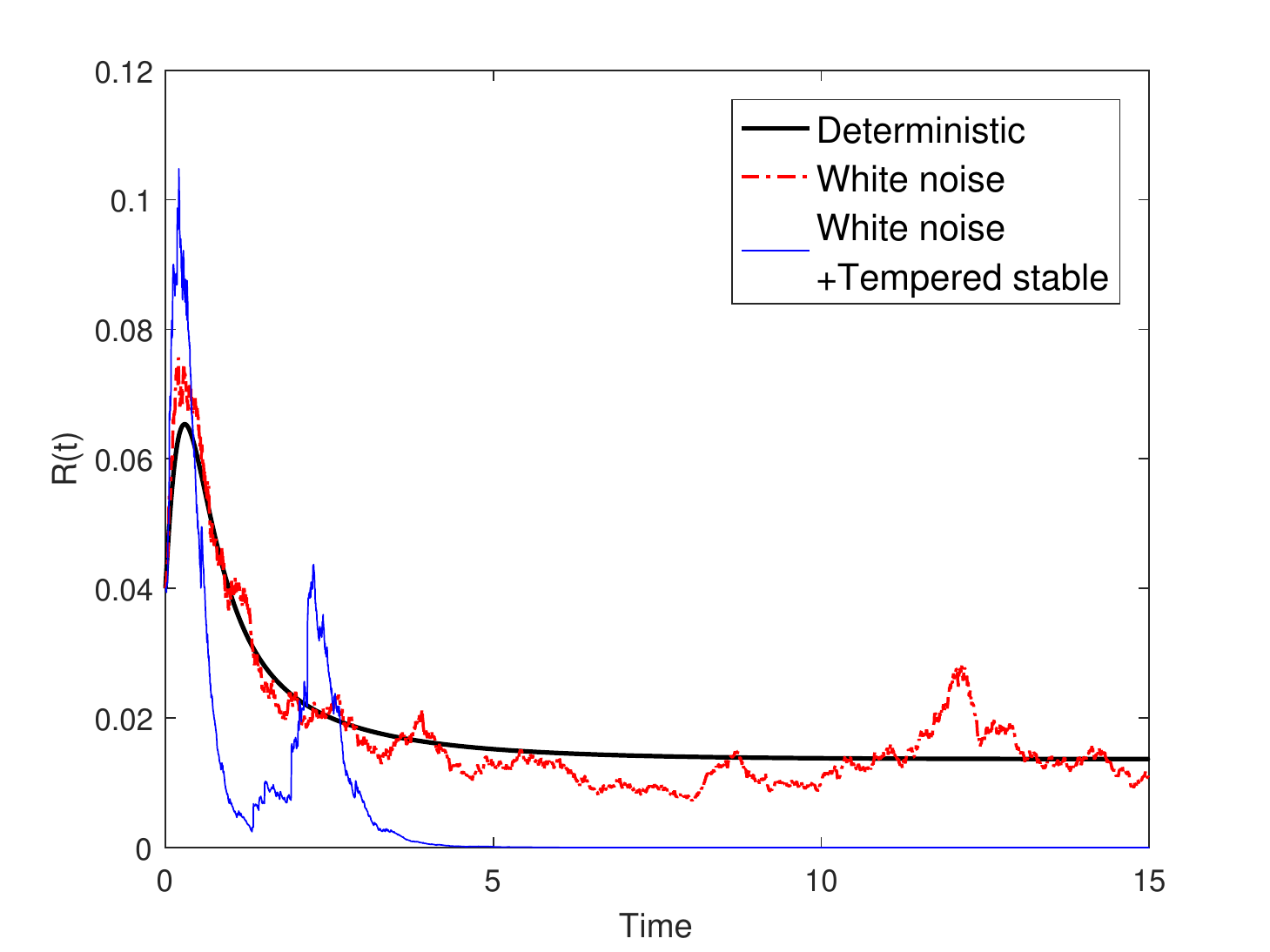}
\caption{\small Extinction of the recovered population} 
\end{subfigure}
\caption{\small Disease extinction in the epidemic population model with $\alpha = 0.7$.} 
\label{fig2}
\end{figure}

\vspace{-0.2cm} 

\begin{figure}[H]
\centering
\hskip-1.4cm
\begin{subfigure}{.5\textwidth}
\centering
\includegraphics[width=1.05\textwidth]{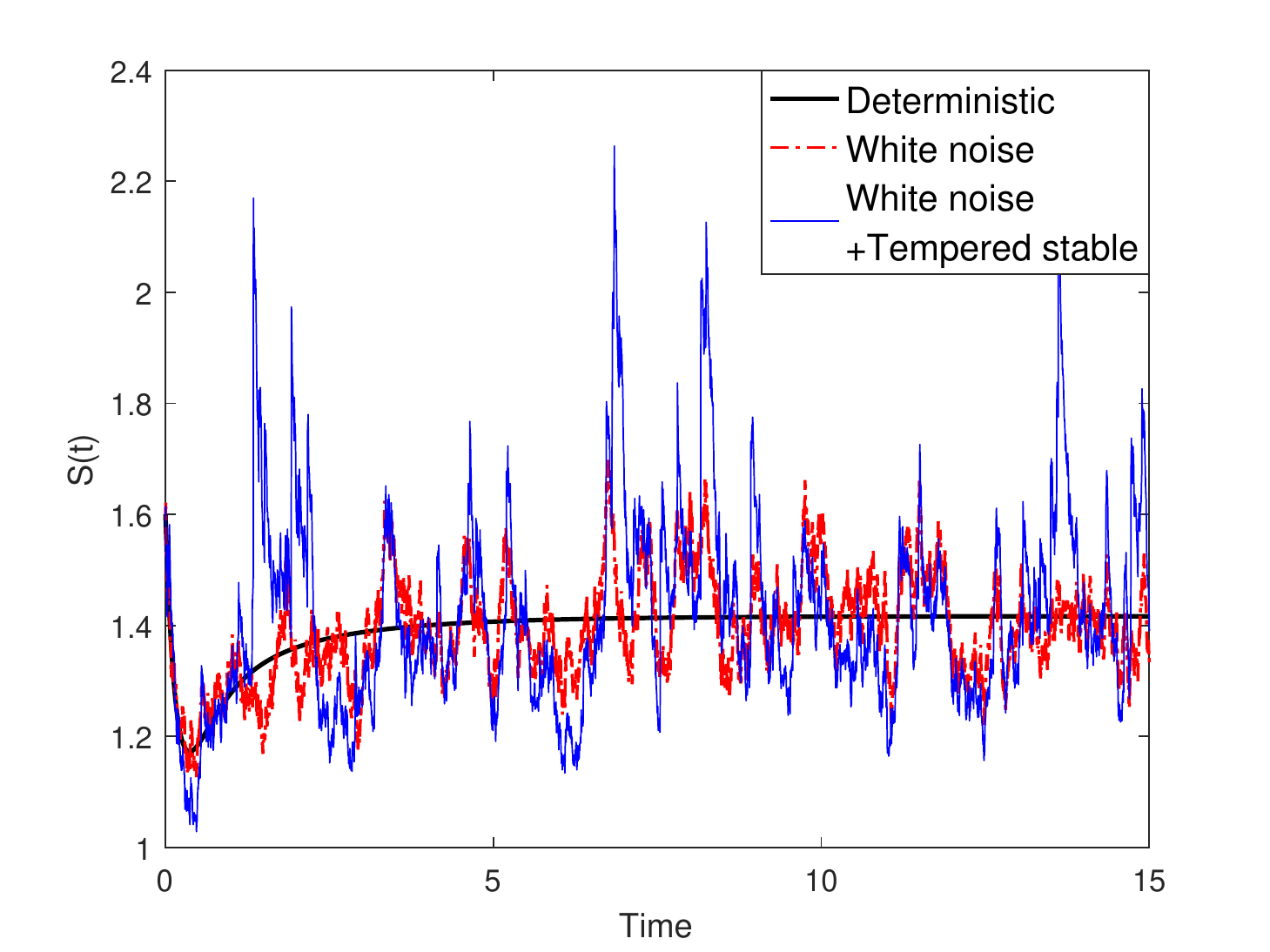}
\caption{\small Dynamical behavior of $S(t)$} 
\end{subfigure}
\hskip-0.5cm
\begin{subfigure}{.5\textwidth}
\centering
\vskip-0.1cm
\includegraphics[width=1.2\textwidth,height=7.05cm]{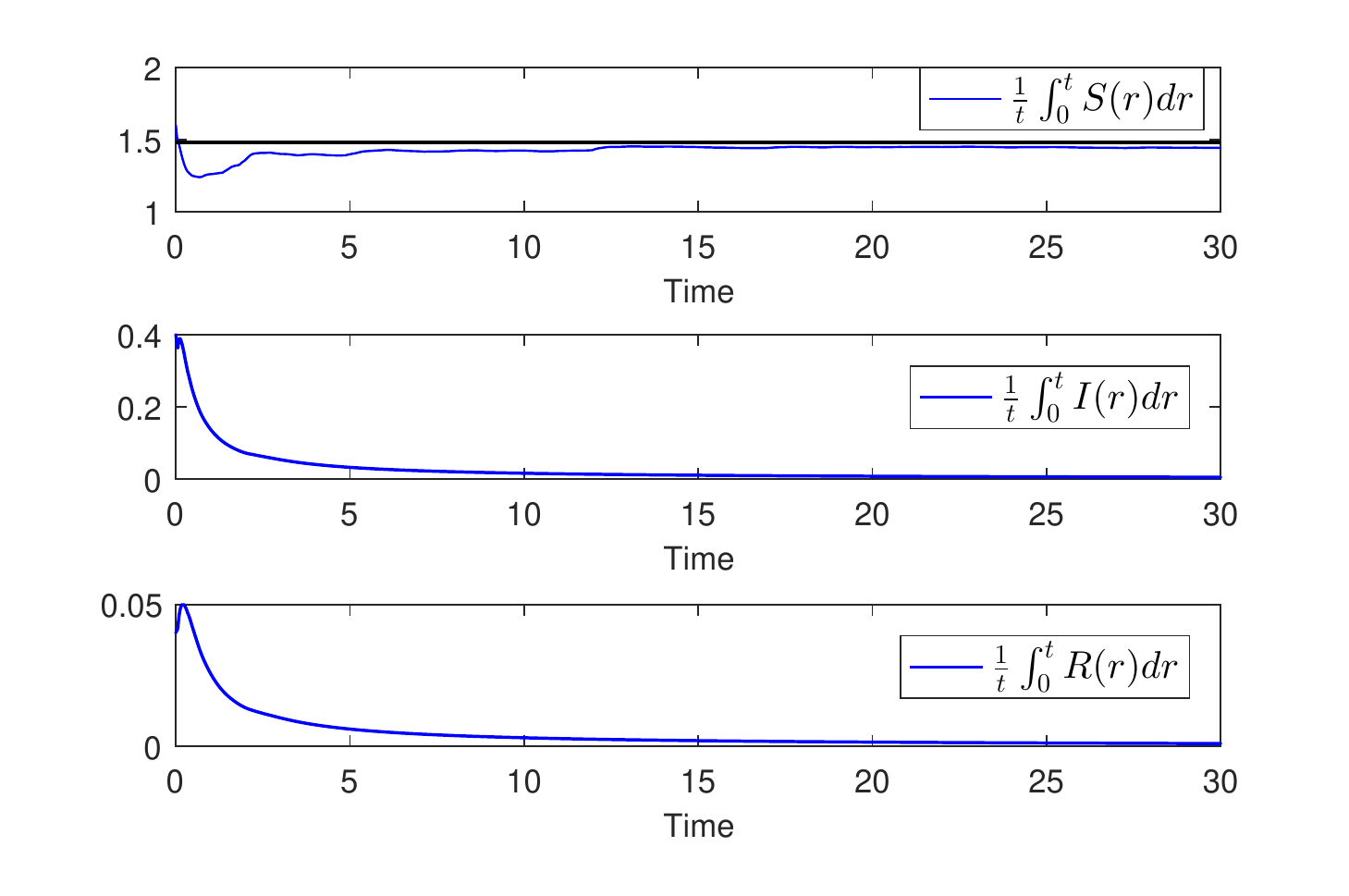}
\vskip-0.5cm
\caption{\small Time averages of $S(t)$, $I(t)$ and $R(t)$} 
\end{subfigure}
\caption{\small Disease extinction in the epidemic population model with $\alpha = 0.7$.} 
\label{fig3}
\end{figure}

\noindent
 We note that the deterministic system is persistent 
  as $\mathcal{R}_0=1.0655>1$, 
with the positive equilibrium value 
$E^\ast=(S^\ast,I^\ast,R^\ast)=(1.417,
0.0723, 
0.0136)$.
On the other hand, for the stochastic system
with $\alpha = 0.7$ we have $\mathcal{\widebar{R}}_0=0.9976<1$, 
and disease extinction is induced by the jump noise with
\begin{equation}\notag
  \lim\limits_{t\rightarrow\infty}\langle S\rangle_t=\frac{\Lambda}{\mu} = 1.5,
  \quad \lim\limits_{t\rightarrow\infty}I_t=0,\quad \lim\limits_{t\rightarrow\infty}R_t=0, \quad \mathbb{P}\mbox{-}a.s. 
\end{equation}
according to Theorem~\ref{t4.1}, see Figures~\ref{fig2}-\ref{fig3}. 

\vspace{-0.4cm}

\vspace{-0.3cm}
 
\noindent 
We also note that the tempered stable model generates jumps of large size
which can model sudden disease outbreak.
Next, we decrease the value of the index to $\alpha=0.2$ and keep the initial value and other parametric values unchanged, in which case Condition~$(H_3^{(2)})$ still holds true
and $\mathcal{\widebar{R}}_0=1.00767>1$.

\begin{figure}[H]
\centering
\hskip-0.2cm
\begin{subfigure}{.5\textwidth}
\centering
\includegraphics[width=1\textwidth]{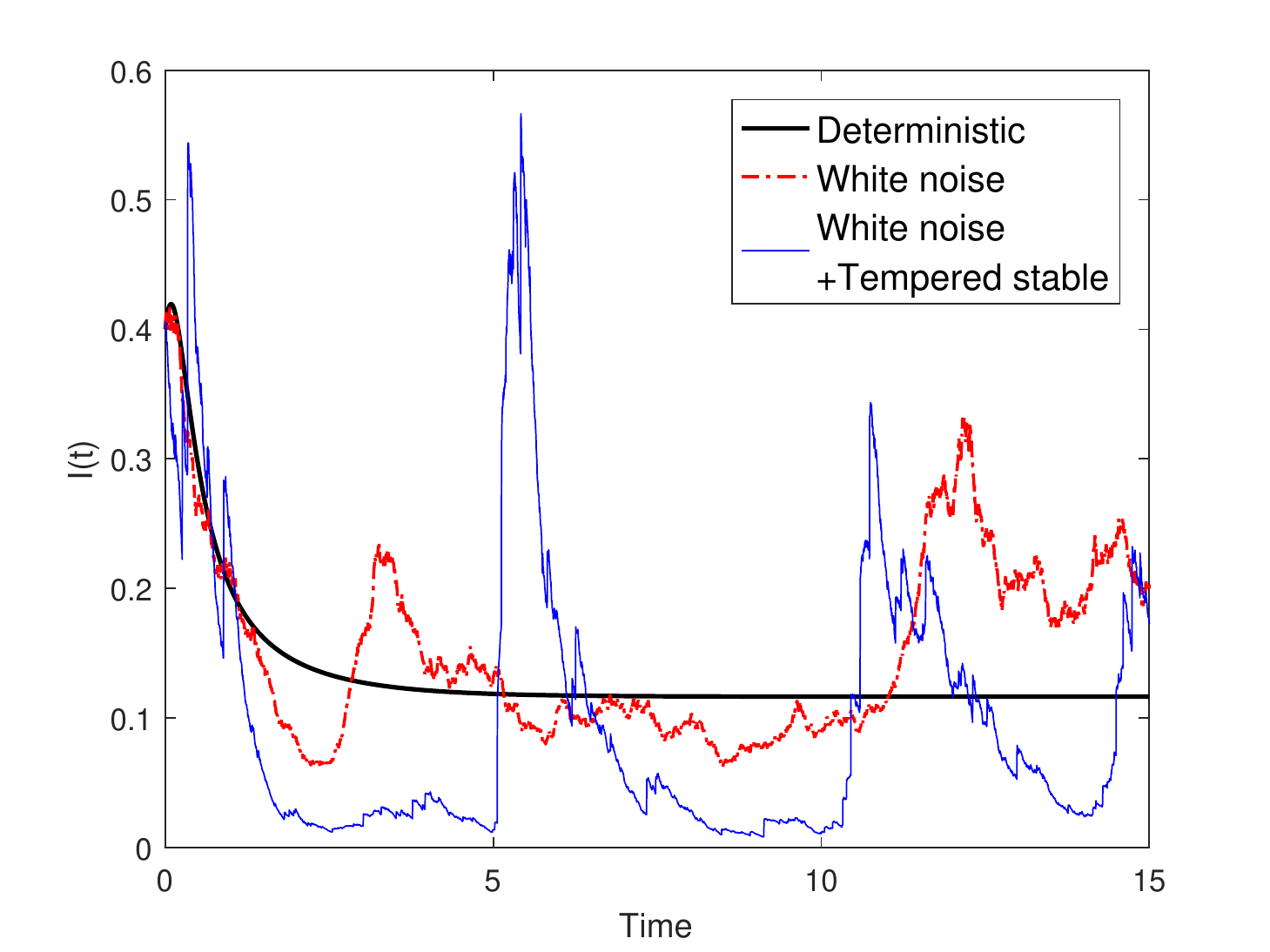}
\caption{\small Persistence of the infected population} 
\end{subfigure}
\begin{subfigure}{.5\textwidth}
\centering
\includegraphics[width=1\textwidth]{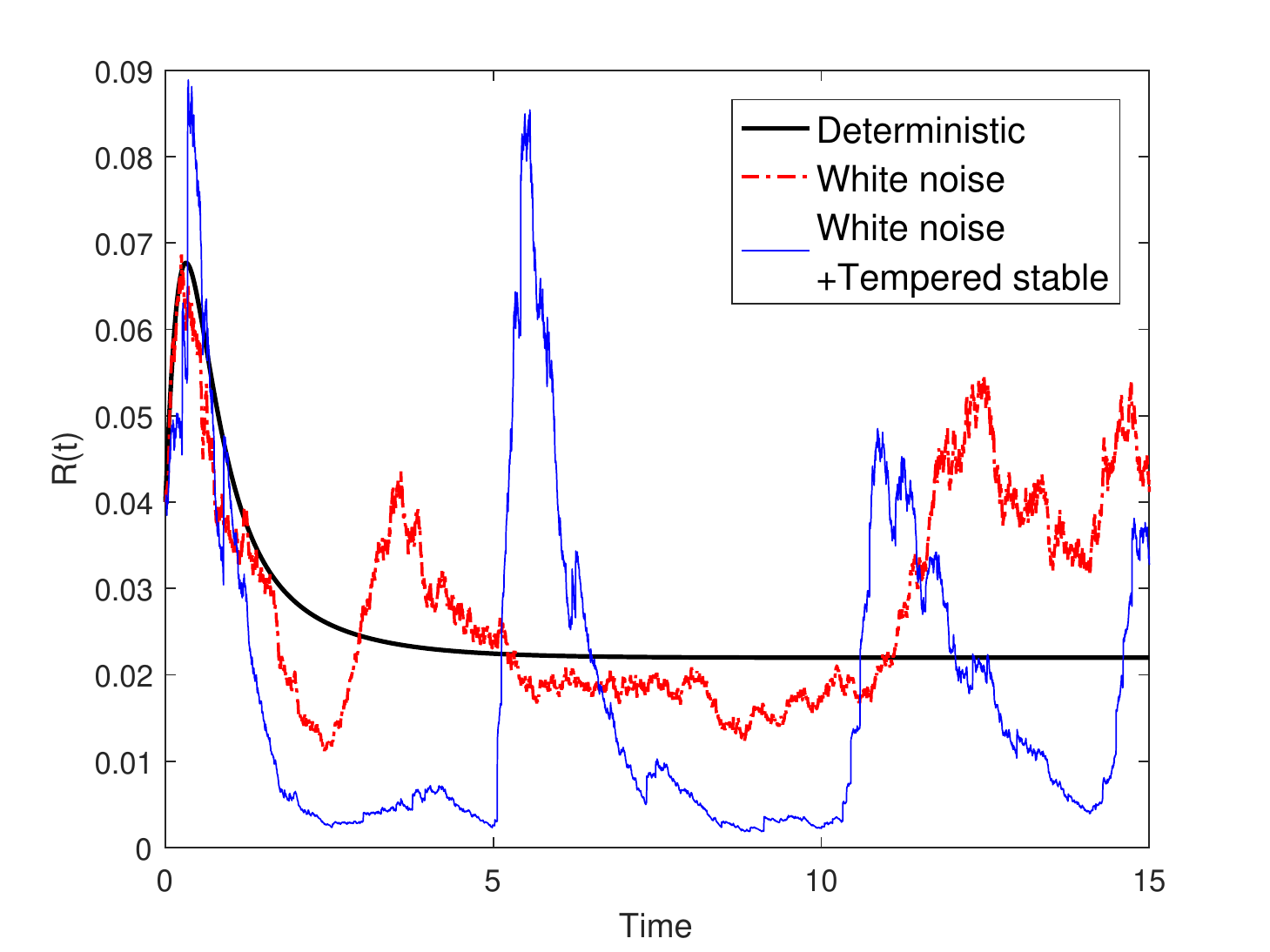}
\caption{\small Persistence of the recovered population} 
\end{subfigure}
\caption{\small Persistence in the epidemic population model with $\alpha = 0.2$.} 
\label{fig4}
\end{figure}

\vspace{-0.2cm} 

\begin{figure}[H]
\centering
\hskip-1.4cm
\begin{subfigure}{.5\textwidth}
\centering
\includegraphics[width=1.05\textwidth]{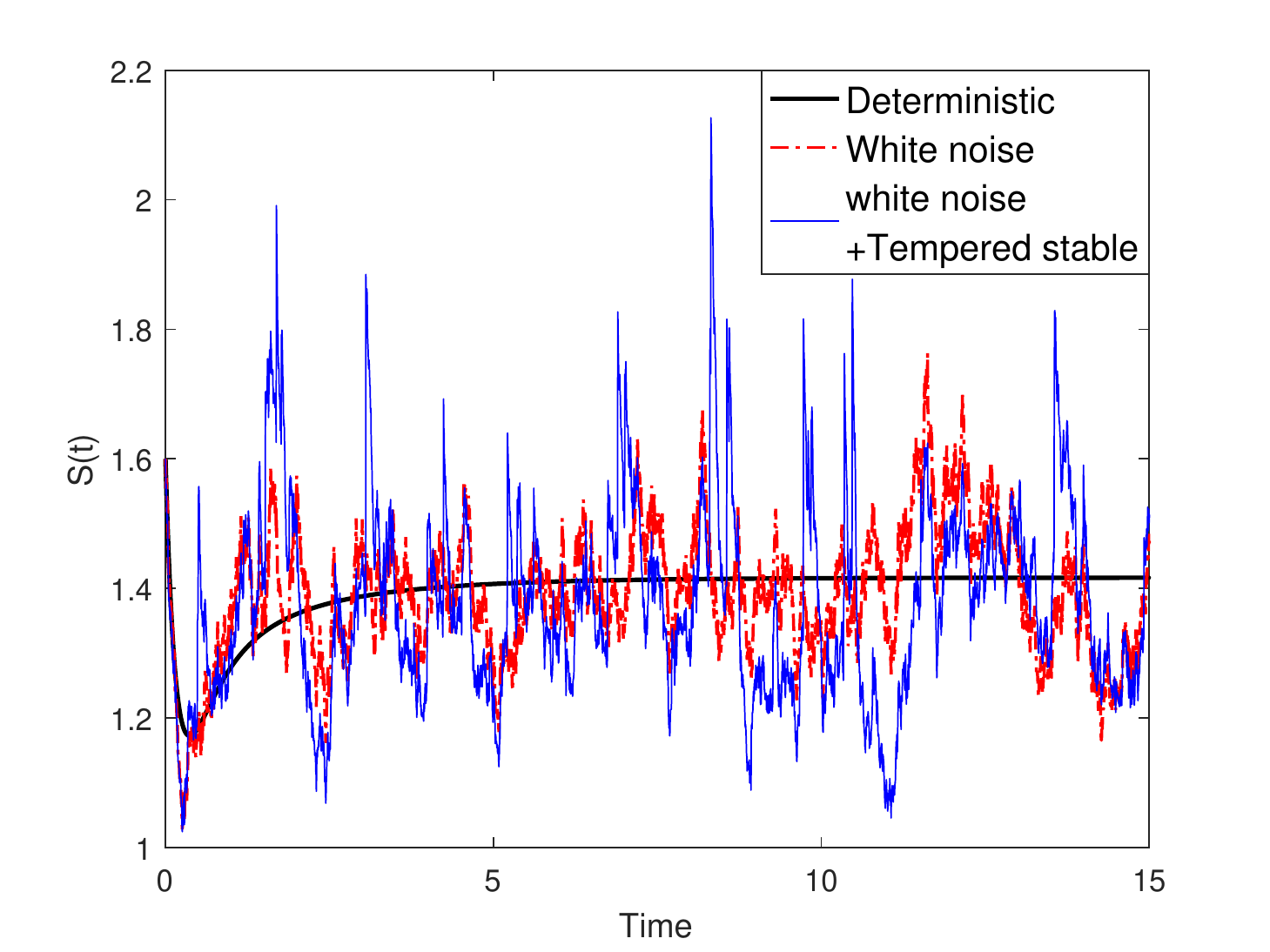}
\caption{\small Dynamical behavior of $S(t)$} 
\end{subfigure}
\hskip-0.5cm
\begin{subfigure}{.5\textwidth}
\centering
\vskip-0.2cm
\includegraphics[width=1.2\textwidth,height=6.8cm]{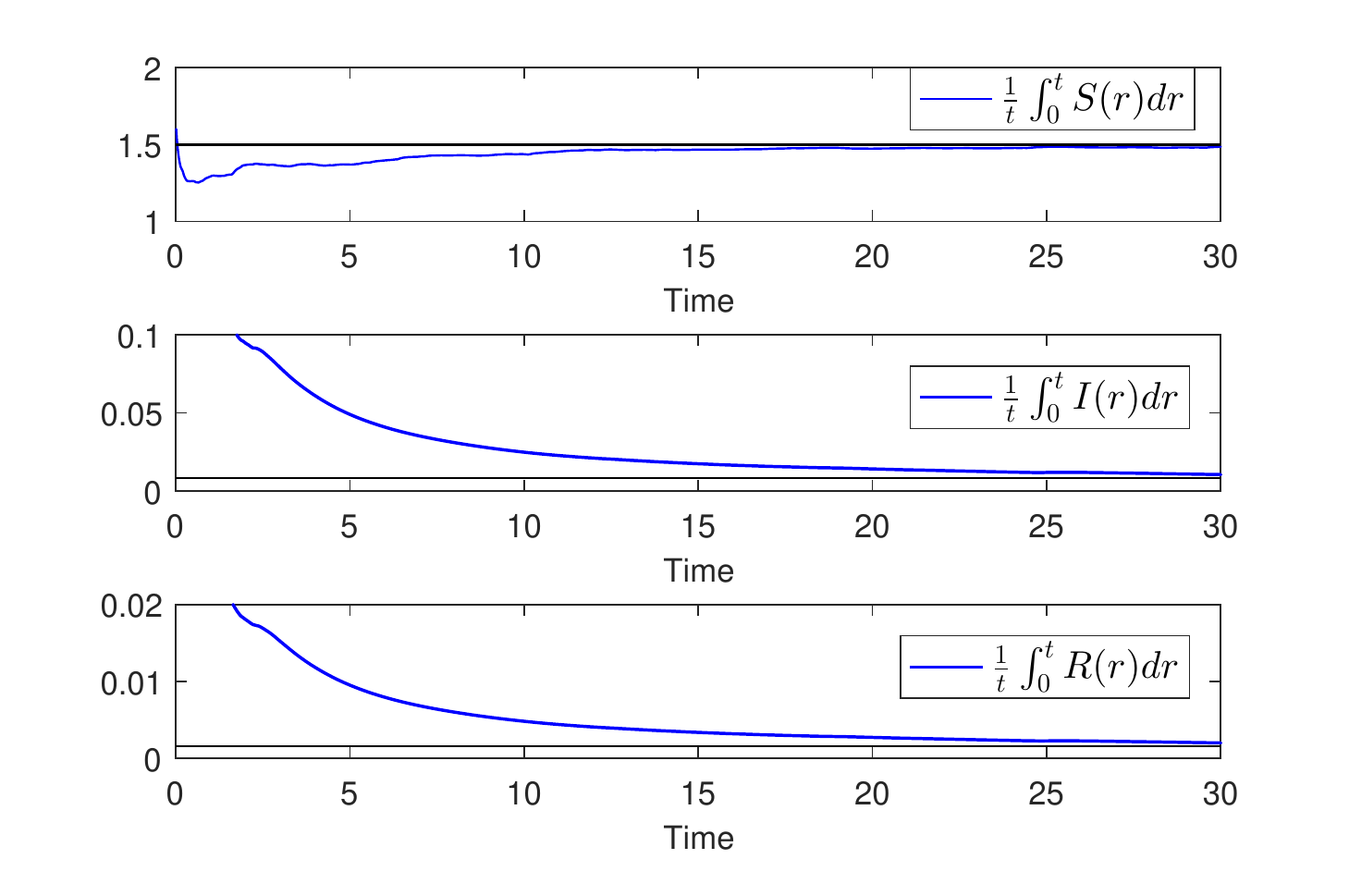}
\vskip-0.5cm
\caption{\small Time averages of $S(t)$, $I(t)$ and $R(t)$} 
\end{subfigure}
\caption{\small Persistence in the epidemic population model with $\alpha = 0.2$.} 
\label{fig5}
\end{figure}

\noindent
 Based on Theorem~\ref{t4.2}, the solution $(S_t,I_t,R_t)$ of
stochastic system (1.1a)-(1.1c) satisfies
$\lim\limits_{t\rightarrow\infty}\langle S\rangle_t=1.49$,
$\lim\limits_{t\rightarrow\infty}\langle I\rangle_t=0.0085$, 
$\lim\limits_{t\rightarrow\infty}\langle R\rangle_t=0.0016$.
The system is persistent and 
 the disease becomes endemic, as illustrated in Figures~\ref{fig4}-\ref{fig5}. 

  \vspace{-0.3cm}

\vspace{-0.2cm}

\noindent
Finally, we consider a pure jump model with two different values
$\alpha^{(1)}$ and $\alpha^{(2)}$ 
and L\'evy measures $\nu^{(1)} (dz)$ and
$\nu^{(2)} (dz)$ given by \eqref{measure3} as 
$$ 
 \nu^{(j)} (dz )
 =
  \frac{ k_+}{z^{\alpha^{(j)} +1}}
 \re^{-\lambda_+ z } dz, \qquad j = 1,2,
$$ 
 while normalizing the jump size variances 
 $$ 
  \big( \sigma^{(1)}_i\big)^2 \int_{\real^3 \setminus\{0\}}
  z^2 {\nu}(dz)
  =
    \big( \sigma^{(2)}_i\big)^2 \int_{\real^3 \setminus\{0\}}
  z^2 {\nu}(dz)
$$ 
  to the same level in both cases, i.e. 
   $$ 
  k_+ \big(\sigma^{(1)}_i\big)^2
\frac{\Gamma(2-\alpha^{(1)})}{\lambda_+^{2-\alpha^{(1)}}}
  =
  k_+ \big(\sigma^{(2)}_i\big)^2
\frac{\Gamma(2-\alpha^{(2)})}{\lambda_+^{2-\alpha^{(2)}}}
$$ 
 with $k_+=2.8$, $\lambda_+=1.2$.
 When $\alpha^{(1)}=0.2$ we take
$\sigma_1^{(1)}=0.2$, $\sigma_2^{(1)}=0.8$ and $\sigma_3^{(1)}=0.5$,
in which case we have $\mathcal{\widebar{R}}_0=1.01>1$
and both $I(t)$ and $R(t)$ are persistent.

\begin{figure}[H]
\centering
\hskip-0.2cm
\begin{subfigure}{.5\textwidth}
\centering
\includegraphics[width=1\textwidth]{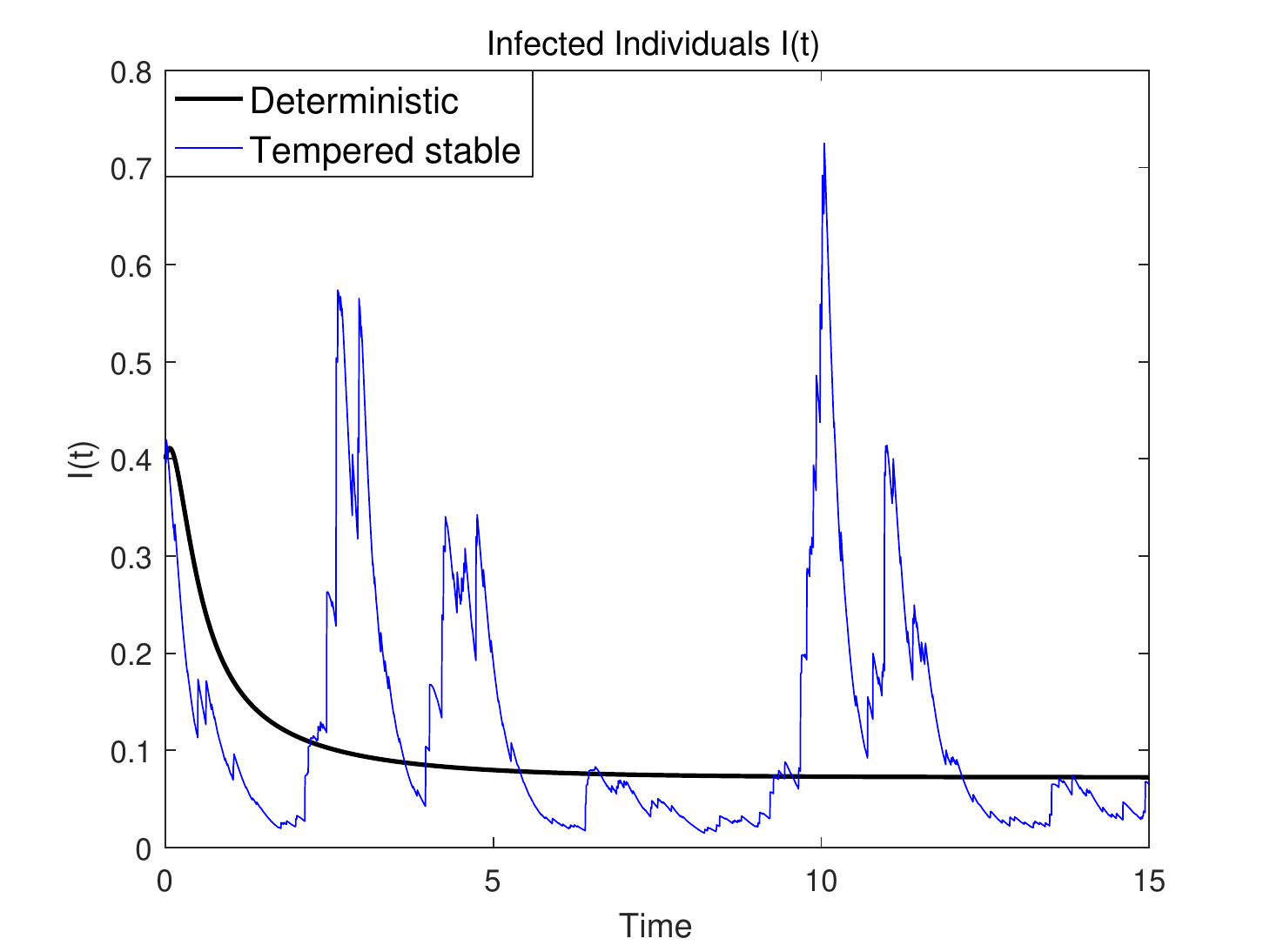}
\caption{\small Persistence of $I(t)$ for $\alpha^{(1)} =0.2$} 
\end{subfigure}
\begin{subfigure}{.5\textwidth}
\centering
\includegraphics[width=1\textwidth]{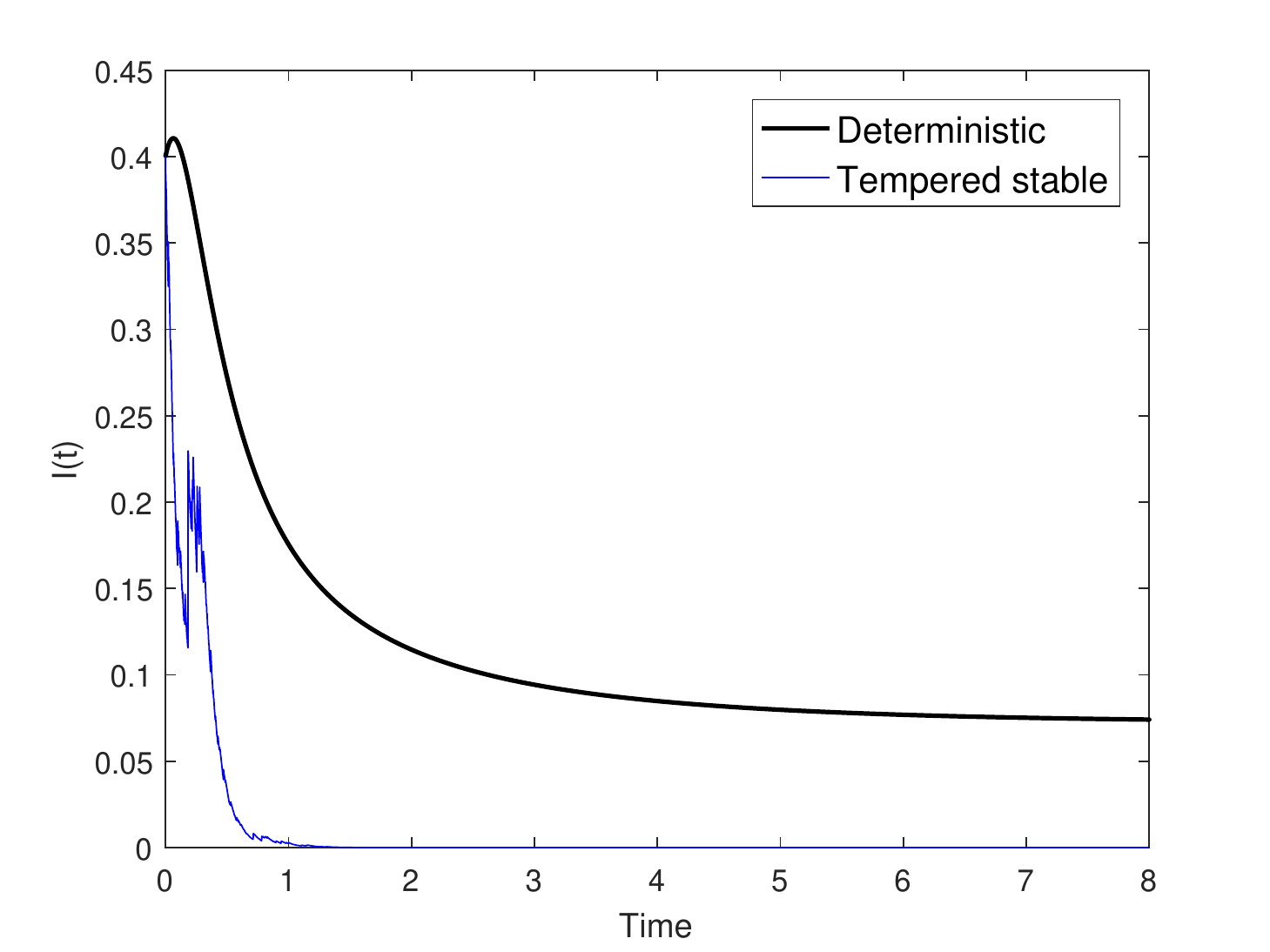}
\caption{\small Extinction of $I(t)$ for $\alpha^{(2)} = 0.9$} 
\end{subfigure}
\caption{\small Behavior of the infected population for two different values of $\alpha$.}\label{fig6}
\end{figure}

\begin{figure}[H]
\centering
\hskip-0.2cm
\begin{subfigure}{.5\textwidth}
\centering
\includegraphics[width=1\textwidth]{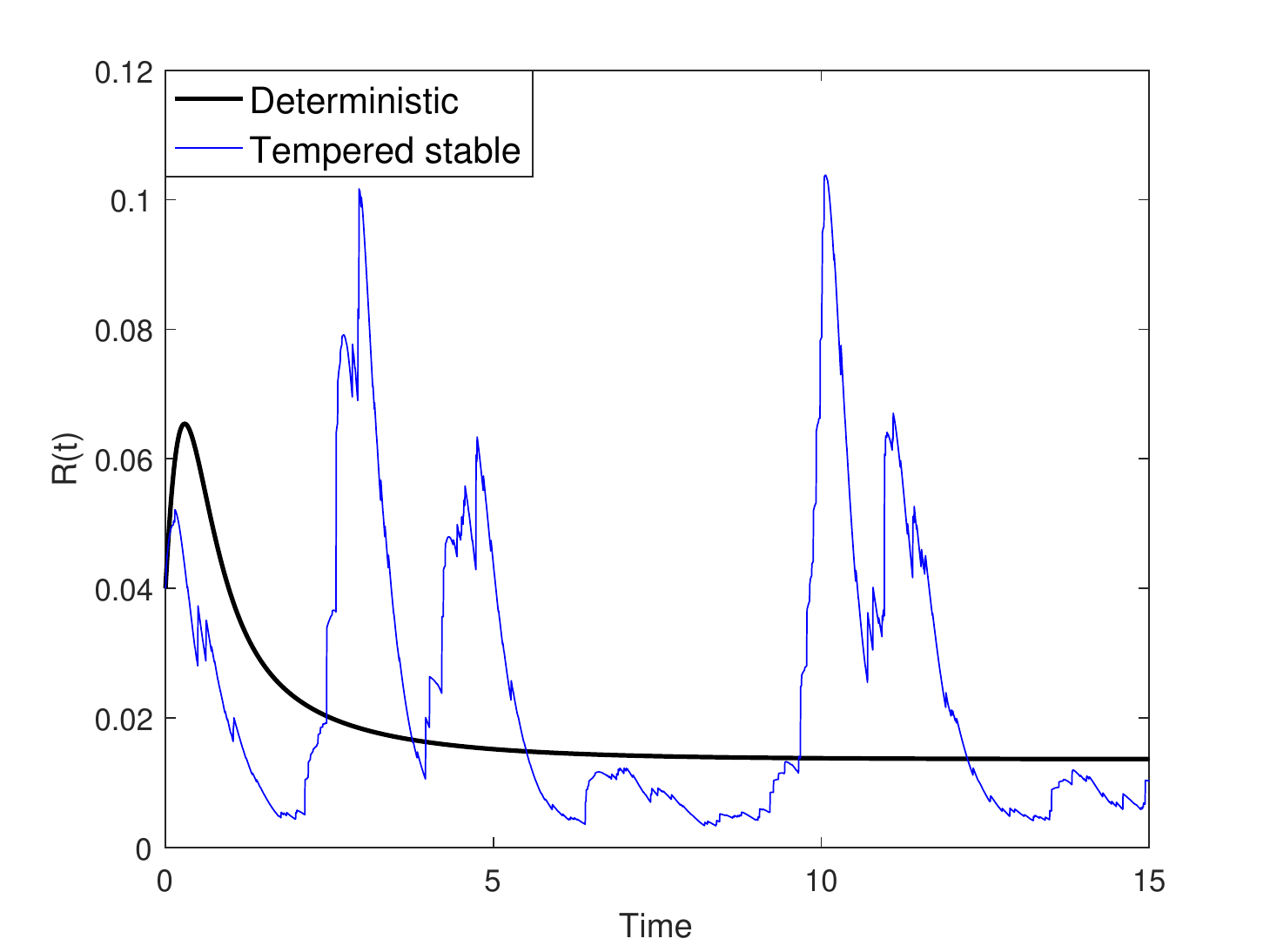}
\caption{\small Persistence of $R(t)$ for $\alpha^{(1)} = 0.2$} 
\end{subfigure}
\begin{subfigure}{.5\textwidth}
\centering
\includegraphics[width=1\textwidth]{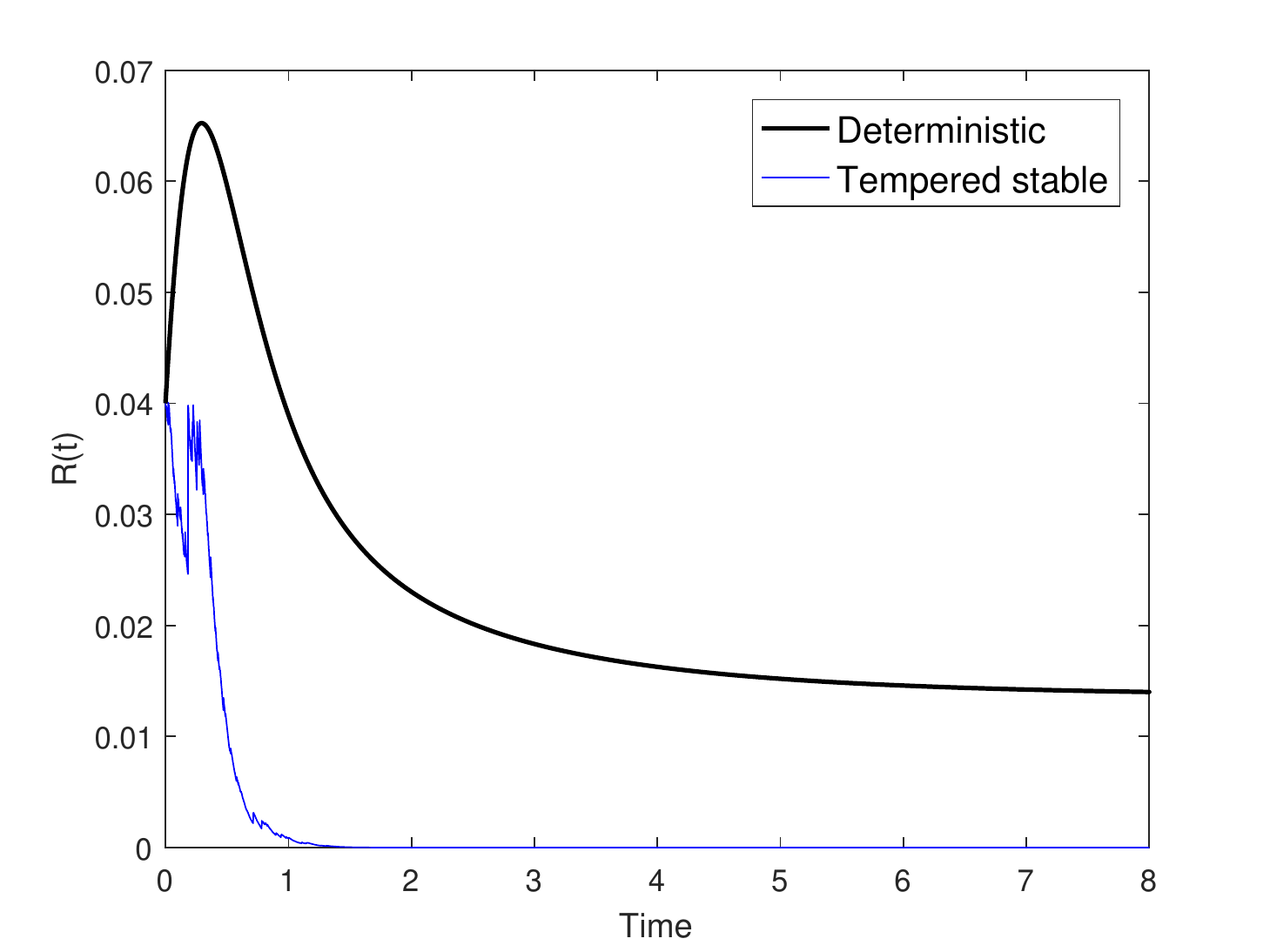}
\caption{\small Extinction of $R(t)$ for $\alpha^{(2)} = 0.7$} 
\end{subfigure}
\caption{\small Behavior of the recovered population for two different values of $\alpha$.}\label{fig7}
\end{figure}

\noindent
 When $\alpha^{(2)}=0.9$ we take
$\sigma_1^{(2)}=0.1857$,
$\sigma_2^{(2)}=0.7426$ and 
$\sigma_3^{(2)}=0.4641$,
in which case we have $\mathcal{\widebar{R}}_0=0.99<1$, 
and both $I(t)$ and $R(t)$ become extinct,
showing that persistence and extinction can depend
on the shape of the jump size distribution for a given variance level,
see Figures~\ref{fig6}-\ref{fig7}. 
In particular, the presence of larger positive jumps for small values
of $\alpha$ can result into persistence of the disease. 

\section{Appendix}
  This section is devoted to proof
  arguments which are similar to the literature, see \cite{wulia}. 
The next Lemma~\ref{Lemma 2.2}, which extends
Lemma~2 of \cite{mazhien} to possibly discontinuous functions $f$,
is needed for the proof of
Theorem~\ref{t4.2}.
\begin{lemma}\label{Lemma 2.2}
  Let $f: \real_+ \to \real_+$ be a function
  integrable on any interval $[0,t]$, $t>0$,
  and consider $\Phi : \real_+ \to \real$ a function such that
  $\lim\limits_{t\rightarrow\infty} ( \Phi(t) / t ) =0$.
\begin{enumerate}[i)]
\item Assume that there exist nonnegative
 constants $\rho_0 \geq 0$, $T\geq0$ such that
\begin{equation}\notag
\log f(t)\leq\rho t-\rho_0\int_0^tf(r)dr+\Phi(t), \quad  \mathbb{P}\mbox{-}a.s.
\end{equation}
for all $t\geq T$, where $\rho \in \real$. Then we have
$$
\left\{
\begin{array}{ll}
  \displaystyle
  \limsup\limits_{t\rightarrow\infty}\frac{1}{t}\int_0^tf(r)dr\leq\frac{\rho}{\rho_0}, \quad  \mathbb{P}\mbox{-}a.s.
  & \mbox{if } \quad \rho\geq0;
  \\
  \\
  \displaystyle
  \lim\limits_{t\rightarrow\infty}f(t)=0, \quad  \mathbb{P}\mbox{-}a.s. & \mbox{if } \quad \rho<0.
\end{array}
\right.
$$
\item Assume that there exists
  positive constants $\rho$, $\rho_0$, and $T\geq0$ such that
\begin{equation}\notag
\log f(t)\geq\rho t-\rho_0\int_0^tf(r)dr+\Phi(t), \quad  \mathbb{P}\mbox{-}a.s.
\end{equation}
for all $t\geq T$. Then we have
\begin{equation}\notag
\liminf\limits_{t\rightarrow\infty}\frac{1}{t}\int_0^tf(r)dr\geq\frac{\rho}{\rho_0}, \quad  \mathbb{P}\mbox{-}a.s.
\end{equation}
\end{enumerate}
\end{lemma}
\begin{Proof}
Define
$$F(t)=\int_0^tf(r)dr.$$
By the integrability of $f(t)$, 
the Lebesgue differentiation theorem (see e.g. Theorem 1.6.11 in \cite{tao})
shows that $F(t)$ is continuous and almost everywhere differentiable,
with
$$\frac{dF(t)}{dt}=f(t)$$
for almost every $t\geq0$.
The rest of the proof is the same as in Lemma 2 of \cite{mazhien}, 
and is omitted.
\end{Proof}
\begin{Proofy}\ref{t4.1}.
  The proof of existence and uniqueness of solutions
  for stochastic differential equation driven by L\'evy processes in \cite{applebk2}
  ensure the integrability of $S_t$, $I_t$ and $R_t$
  on any bounded interval $[0,T]$.
  In view of \eqref{lsir1}-\eqref{lsir3}, we deduce
\begin{eqnarray*}
  \lefteqn{
    \! \! \! \! \! \! \! \! \! \! \! \! \!
    \! \! \! \! \!
    \frac{S_t-S_0}{t}+\frac{I_t-I_0}{t}
= \Lambda-\mu\langle S\rangle_t-(\mu+\varepsilon+\eta)\langle I\rangle_t+\frac{1}{t}\int_0^tS_rdB^\varrho_1(r)+\frac{1}{t}\int_0^tI_rdB^\varrho_2(r)
  }
  \\
& &+\frac{1}{t}\int_0^tS_{r^-}\int_{\real^3 \setminus\{0\}}\gamma_1(z)\tilde{N}(dr,dz)
+\frac{1}{t}\int_0^tI_{r^-}\int_{\real^3 \setminus\{0\}}\gamma_2(z)\tilde{N}(dr,dz)
,
\end{eqnarray*}
which yields
\begin{equation}\label{relation}
  \mu\langle S\rangle_t+(\mu+\varepsilon+\eta)\langle I\rangle_t=\Lambda-\varphi (t),
\end{equation}
where
\begin{eqnarray*}
\varphi (t) & :=& \frac{S_t-S_0}{t}+\frac{I_t-I_0}{t}-\frac{1}{t}\int_0^tS_rdB^\varrho_1(r)-\frac{1}{t}\int_0^tI_rdB^\varrho_2(r)\\
& &-\frac{1}{t}\int_0^tS_{r^-}\int_{\real^3 \setminus\{0\}}\gamma_1(z)\tilde{N}(dr,dz)
-\frac{1}{t}\int_0^tI_{r^-}\int_{\real^3 \setminus\{0\}}\gamma_2(z)\tilde{N}(dr,dz).
\end{eqnarray*}
By the It\^{o} formula for L\'evy-type stochastic integrals
(see Theorem~1.16 in \cite{yiteng}) and \eqref{relation},
 letting
 \begin{equation}
   \label{m2t} 
M_2(t):=\int_0^t\int_{\real^3 \setminus\{0\}}\log (1+\gamma_2(z))\tilde{N}(ds,dz)
\end{equation}
 we have
\begin{eqnarray}
\nonumber
  \log I_t & = & \log I_0+\beta\int_0^tS_rdr-(\mu+\varepsilon+\eta)t-\beta_2t
 +B^\varrho_2(t)+M_2(t)
  \\
\nonumber
  &=&\log I_0+\left(
 \beta\frac{\Lambda}{\mu}-(\mu+\varepsilon+\eta + \beta_2 ) \right) t
 -\frac{\beta(\mu+\varepsilon+\eta)}{\mu}\int_0^tI_rdr
 \\
\label{lnI}
 & &
 -\frac{\beta}{\mu}t\varphi (t)
+B^\varrho_2(t)
 +
 M_2(t)
 \\
\nonumber
&\leq&\log I_0+
(\mu+\varepsilon+\eta)(\mathcal{\widebar{R}}_0-1) t -\frac{\beta}{\mu}t\varphi (t)
+B^\varrho_2(t)
+M_2(t).
\end{eqnarray}
 We deduce from Lemmas~\ref{l3.1}, \ref{l3.3} and \ref{l3.4} that
 \begin{equation}
   \label{fai}
  \lim\limits_{t\rightarrow\infty}\varphi (t)=0, \quad  \mathbb{P}\mbox{-}a.s.
\end{equation}
 In addition, under $(H_5)$ we have
 \begin{eqnarray*}
  \int_0^t\frac{d\langle M_2,M_2\rangle(r)}{(1+r)^2}dr
  & = &
  \int_0^t\frac{1}{(1+r)^2}
  dr
  \int_{\real^3 \setminus\{0\}}
  \big( \log (1+\gamma_2(z)) \big)^2
  \nu(dz)
  \\
   & = &
  \frac{t}{1+t}\int_{\real^3 \setminus\{0\}}
  \big(
  \log (1+\gamma_2(z)) \big)^2\nu(dz)<+\infty,
 \quad t\in \real_+,
\end{eqnarray*}
 hence by the law of large numbers for local martingales
(see Theorem~1 in \cite{dashudinglvlevy}) we have
\begin{equation}\label{yang}
  \lim\limits_{t\rightarrow\infty}\frac{M_2(t)}{t}=0, \quad  \mathbb{P}\mbox{-}a.s.
\end{equation}
By the law of large numbers
(see Theorem~3.4 in Chapter~1 of \cite{mao2008}) we also get
\begin{equation}
\label{**}
\lim\limits_{t\rightarrow\infty}\frac{B^\varrho_2(t)}{t}=0, \quad
 \mathbb{P}\mbox{-}a.s.
\end{equation}
 Therefore, by \eqref{lnI}, if $\mathcal{\widebar{R}}_0<1$ we have
\begin{equation}\nonumber 
  \limsup\limits_{t\rightarrow\infty}\frac{\log I_t}{t}\leq (\mu+\varepsilon+\eta)(\mathcal{\widebar{R}}_0-1)<0, \quad  \mathbb{P}\mbox{-}a.s.,
\end{equation}
which, together with the positivity of $I_t$, implies
\begin{equation}\label{ID}
  \lim\limits_{t\rightarrow\infty}I_t=0, \quad  \mathbb{P}\mbox{-}a.s.
\end{equation}
In other words, the disease goes to extinction with probability one. Furthermore,
from \eqref{relation} we obtain
\begin{equation}\notag
  \lim\limits_{t\rightarrow\infty}\langle S\rangle_t=\frac{\Lambda}{\mu}, \quad  \mathbb{P}\mbox{-}a.s.
\end{equation}
 We derive from \eqref{lsir3} that
\begin{equation}\notag
\frac{R_t-R_0}{t}=-\frac{\mu}{t}\int_0^tR_rdr+\frac{\eta}{t}\int_0^tI_rdr+
\frac{1}{t}\int_0^tR_rdB^\varrho_3(r)+\frac{1}{t}\int_0^tR_{r^-}\int_{\real^3 \setminus\{0\}}\gamma_3(z)\tilde{N}(dr,dz),
\end{equation}
 and taking limits on both sides yields
\begin{eqnarray}\label{limr}
\nonumber
\mu \lim\limits_{t\rightarrow\infty}\frac{1}{t}\int_0^tR_rdr
&=&\eta \lim\limits_{t\rightarrow\infty}\frac{1}{t}\int_0^tI_rdr-\lim\limits_{t\rightarrow\infty}\frac{R_t}{t}
+\lim\limits_{t\rightarrow\infty}\frac{1}{t}\int_0^tR_rdB^\varrho_3(r)\\
  & & +\lim\limits_{t\rightarrow\infty}\frac{1}{t}\int_0^tR_{r^-}\int_{\real^3 \setminus\{0\}}\gamma_3(z)\tilde{N}(dr,dz), \quad  \mathbb{P}\mbox{-}a.s.
\end{eqnarray}
Together with \eqref{ID} and the conclusions in Lemmas~\ref{l3.1},
\ref{l3.3} and \ref{l3.4}, we conclude to
\begin{equation}\notag
\lim\limits_{t\rightarrow\infty}R_t=0, \quad  \mathbb{P}\mbox{-}a.s.\end{equation}
\end{Proofy}
\begin{Proofy}\ref{t4.2}.
  By \eqref{relation} we deduce that
\begin{equation}\notag
  \beta\langle S\rangle_t=\frac{\beta\Lambda}{\mu}-\frac{\beta}{\mu}\varphi (t)-\frac{\beta(\mu+\varepsilon+\eta)}{\mu}\langle I\rangle_t.
\end{equation}
It then follows from \eqref{lnI} that
\begin{eqnarray*}
\log I_t&=&\log I_0+\beta\int_0^tS_rdr-(\mu+\varepsilon+\eta + \beta_2 ) t
  \\
  & &+B^\varrho_2(t)+\int_0^t\int_{\real^3 \setminus\{0\}}\log (1+\gamma_2(z))\tilde{N}(ds,dz)\\
& = & \left(
  \beta\frac{\Lambda}{\mu}-(\mu+\varepsilon+\eta)-\beta_2 \right) t
 -\frac{\beta(\mu+\varepsilon+\eta)}{\mu}\int_0^tI_rdr
 + \Psi (t),
\end{eqnarray*}
where we denote
$$
\Psi(t):=\log I_0-\frac{\beta}{\mu}t\varphi (t)+B^\varrho_2(t)+M_2(t),
\qquad t\in \real_+, 
$$
and $M_2(t)$ is defined as in \eqref{m2t}.  
 From \eqref{fai}, \eqref{yang} and \eqref{**}
 it follows that
 $\lim\limits_{t\rightarrow\infty}(\Psi(t)/t)=0$,
 hence
 applying Lemma \ref{Lemma 2.2} to the function $f(t):=I_t$
 which is a.s. integrable over $[0,T]$, $T>0$,
 we obtain
\begin{equation}\nonumber 
  \lim\limits_{t\rightarrow\infty}\langle I\rangle_t=\frac{\mu\left(\beta \Lambda / \mu
    -(\mu+\varepsilon+\eta)-\beta_2 \right)}{\beta(\mu+\varepsilon+\eta)}=\frac{\mu}{\beta}(\mathcal{\widebar{R}}_0-1), \quad  \mathbb{P}\mbox{-}a.s.
\end{equation}
 Consequently, on account of \eqref{relation} and \eqref{fai} we get
\begin{equation}\nonumber 
  \lim\limits_{t\rightarrow\infty}\langle S\rangle_t=\frac{\Lambda}{\mu}-\frac{(\mu+\varepsilon+\eta)}{\beta} (\mathcal{\widebar{R}}_0-1) =S^\ast+\frac{\beta_2}{\beta}, \quad  \mathbb{P}\mbox{-}a.s.,
\end{equation}
and it follows from \eqref{limr}
and Lemmas~\ref{l3.1}, \ref{l3.3} and \ref{l3.4}
that
\begin{equation}\nonumber 
  \lim\limits_{t\rightarrow\infty}\langle R\rangle_t=\frac{\eta}{\beta}(\mathcal{\widebar{R}}_0-1), \quad  \mathbb{P}\mbox{-}a.s.
\end{equation}
\end{Proofy}

\subsubsection*{Conclusion}
In this paper, we consider a stochastic version of the SIR epidemic model
\eqref{lsir1}-\eqref{lsir3}, 
driven by correlated Brownian and L\'{e}vy jump components with 
 heavy tailed increments. 
We present new solution estimates using the parameter $\lambda(p)$ defined in
\eqref{lmbda2} and Kunita's inequality for jump processes in the
key Lemmas~\ref{l3.1} and \ref{l3.4}. 
Our approach relaxes
the restriction on the finiteness of the L\'{e}vy measure $\nu(dz)$
imposed in \cite{amllevy} and \cite{wulia},
and our definition of the parameter $\lambda(p)$ in \eqref{lmbda2}
applies to a wider range of L\'{e}vy measures. 
In Theorems~\ref{t4.1} and \ref{t4.2} we derive
the basic reproduction number $\mathcal{\widebar{R}}_0$ 
which characterizes the extinction and persistence properties
of the stochastic epidemic system
\eqref{lsir1}-\eqref{lsir3}. 
As an illustration we present numerical simulations based on
tempered stable processes, showing that the additional presence of jumps and
the level of the index $\alpha \in (0,1)$ can have a significant influence on the
dynamical behavior of the epidemic system.
 
\footnotesize 

\setcitestyle{numbers}

\def\cprime{$'$} \def\polhk#1{\setbox0=\hbox{#1}{\ooalign{\hidewidth
  \lower1.5ex\hbox{`}\hidewidth\crcr\unhbox0}}}
  \def\polhk#1{\setbox0=\hbox{#1}{\ooalign{\hidewidth
  \lower1.5ex\hbox{`}\hidewidth\crcr\unhbox0}}} \def\cprime{$'$}

\end{document}